\documentclass[envcountsect,envcountsame,oribibl,orivec]{llncs}
\usepackage{amssymb}
\usepackage{amsmath}
\usepackage{url}
\usepackage{enumitem}
\usepackage{colonequals}

\usepackage{graphics}
\usepackage{graphicx}
\usepackage{hyperref}

\usepackage{qtree}

\usepackage{bussproofs}
\def\fCenter{\hspace{1pt}\Rightarrow\hspace{1pt}}
\EnableBpAbbreviations

\usepackage{etoolbox}
\makeatletter
\let\llncs@addcontentsline\addcontentsline
\patchcmd{\maketitle}{\addcontentsline}{\llncs@addcontentsline}{}{}
\patchcmd{\maketitle}{\addcontentsline}{\llncs@addcontentsline}{}{}
\patchcmd{\maketitle}{\addcontentsline}{\llncs@addcontentsline}{}{}
\makeatother
\usepackage{hyperref}
\usepackage{bookmark}
\makeatletter
\newcommand{\todo}[1]{\marginpar{\textbf{TODO\footnotemark}}\@latex@warning{TODO: #1}\footnotetext{ #1}}
\makeatother
\usepackage{tikz}

\usepackage{tikz}
\usetikzlibrary{trees,arrows,automata,shapes,snakes,decorations,decorations.pathmorphing}

\newenvironment{mytikz}[1][0em]{
	\begin{tikzpicture}[>=latex,auto,node distance=6.75em,baseline={([yshift=#1]current bounding box.east)}]
	
	\tikzstyle{w}=[draw,circle,minimum size=3em]
	
	\tikzstyle{e}=[draw,minimum size=2.5em,node distance=3cm]
	
	\tikzstyle{d}=[densely dashed]
	
	\tikzstyle{every edge}=[draw,font=\footnotesize]
	
	\tikzstyle{every label}=[font=\footnotesize]
	
	\tikzstyle{zz}=[decorate,decoration={snake,post=lineto,post length=6pt}]
	
	\tikzstyle{ev}=[anchor=west,node distance=2em]
	
	\tikzstyle{bigger}=[minimum size=4em]
	
	\tikzstyle{bit-farther}=[node distance=2.3cm]
	
	\tikzstyle{farther}=[node distance=3cm]
	
	\tikzstyle{farthest}=[node distance=4cm]
	
	\tikzstyle{l}=[node distance=3.5em]
}{\end{tikzpicture}}

\newcommand{\Prop}{\textsf{Prop}}
\newcommand{\Formulae}{\textsf{Fm}}

\newcommand{\lnegRelevant}{\sim}
\newcommand{\lfussion}{\circ}

\newcommand{\justVarOne}{\mathsf{x}}
\newcommand{\justVarTwo}{\mathsf{y}}

\newcommand{\justConsThree}{\mathsf{c}}

\newcommand{\lrule}[2]{\displaystyle{\frac{#1}{#2}}}
\newcommand{\mprule}{\ensuremath{(\textsf{MP})}}

\newcommand{\RCEArule}{\ensuremath{(\textsf{RCEA})}}
\newcommand{\RCKrule}{\ensuremath{(\textsf{RCK})}}
\newcommand{\RCNrule}{\ensuremath{(\textsf{RCN})}}

\newcommand{\limplies}{\supset}
\newcommand{\limpliesCounterfactual}{>}
\newcommand{\limpliesRelevant}{\rightarrow}
\newcommand{\limpliesCounterRel}{\rightsquigarrow}

\newcommand{\liff}{\equiv}
\newcommand{\powerset}{\wp}

\newcommand{\Logic}[1]{\mathsf{#1}} 

\newcommand{\LP}{\textsf{LP}}
\newcommand{\LPC}{\textsf{LPC}^+}
\newcommand{\LPCint}{\textsf{LPC}^{int}}
\newcommand{\LPCstandard}{\textsf{LPC}^{\prime}}
\newcommand{\JRC}{\textsf{JRC}}
\newcommand{\LPCK}{\textsf{LPCK}^+}
\newcommand{\JFourC}{\textsf{J4C}^+}
\newcommand{\JC}{\textsf{JC}^+}
\newcommand{\JTC}{\textsf{JTC}^+}

\newcommand{\CTerms}{\textsf{Const}}
\newcommand{\VTerms}{\textsf{Var}}
\newcommand{\Terms}{\textsf{Tm}}
\newcommand{\TermsJCR}{\textsf{Term}}
\newcommand{\FormulaeJCR}{\textsf{Form}}

\newcommand{\jbox}[1]{\left[#1\right]\!}

\newcommand{\truthsetModel}[1]{\left\|#1\right\|_\model\!}
\newcommand{\truthset}[2]{\left\|#1\right\|_#2\!}

\newcommand{\tapp}{\cdot}
\newcommand{\ttsum}{+}
\newcommand{\tinspect}{!}

\newcommand{\prop}[1]{Prop(#1)\!}
\newcommand{\Sub}[1]{Sub(#1)\!}
\newcommand{\lset}[1]{\{ #1 \}}

\newcommand{\lknows}{\Box}
\newcommand{\JTB}{\mathsf{JTB}}

\newcommand{\CS}{\textsf{CS}}

\newcommand{\worlds}{W}
\newcommand{\worldsNormal}{W_N}
\newcommand{\worldsNonNormal}{\worlds \setminus \worldsNormal}

\newcommand{\valuation}{V}
\newcommand{\valuationR}{\mathcal{V}}
\newcommand{\entails}{\vDash}

\newcommand{\model}{\mathcal{M}}
\newcommand{\accrel}{R}
\newcommand{\relRelevant}{R}
\newcommand{\accrelTrem}{R_{Tm}}
\newcommand{\accrelFormula}{R_{Fm}}
\newcommand{\relTrem}{R_t}
\newcommand{\relFormula}{R_\phi}

\newcommand{\relModel}[1]{R_{#1}\,}
\newcommand{\relTableau}[1]{r_{#1}\,}

\newcommand{\N}{\mathbb{N}}
\newcounter{enumsave}
\renewcommand{\phi}{\varphi}

\newcommand{\MCS}{\mathsf{MCS}}

\begin{document}

\title{Justification logics with counterfactual and relevant conditionals }
\author{Meghdad Ghari\thanks{This research was in part supported by a grant from IPM (No. 1400030416) and carried out in IPM-Isfahan Branch.}}
\institute{Department of Philosophy, Faculty of Literature and Humanities,\\
	University of Isfahan, Isfahan, Iran \\ and \\ School of Mathematics,
	Institute for Research in Fundamental Sciences (IPM), \\ P.O.Box: 19395-5746, Tehran, Iran \\ \email{m.ghari@ltr.ui.ac.ir}
}

\maketitle

\begin{abstract}

The purpose of this paper is to introduce justification logics based on conditional logics. We introduce a new family of logics, called conditional justification logics, which incorporates a counterfactual conditional in its language. For the semantics, we offer relational models that merge the selection-function semantics of conditional logics with the relational semantics of justification logics. As an application, we formalize Nozick's counterfactual conditions in his analysis of knowledge and investigate their connection to Aumann's concepts of knowledge. Additionally, we explore Gettier's counterexamples to the justified true belief analysis, as well as McGinn's counterexamples to Nozick's analysis of knowledge. Furthermore, we introduce a justification logic that includes a relevant counterfactual conditional and we demonstrate the variable-sharing property for this conditional. We also develop a tableau system for this logic and establish its completeness with respect to Routley relational models. Finally, we formalize Nozick's counterfactual conditions using this relevant counterfactual conditional and represent the sheep in the field example of Chisholm within this logic.

\end{abstract}

\section{Introduction}
\label{sec:Introduction}

Justification logics are a family of epistemic logics that provide justifications and evidence for an agent's belief. The first such logic was Artemov's logic of proofs \cite{Art95TR,Art01BSL}. Since then,   various justification logics have been introduced   in the literature (see e.g. \cite{Artemov-Fitting-Studer-SEP,Art-Fit-Book-2019,Kuz-Stu-Book-2019}).  These logics extend classical propositional logic with  expressions of the form $t:\phi$, where $\phi$ is a formula and $t$ is a justification term, interpreted  as ``$t$ is a justification for $\phi$".  Most justification logics are based on classical logic, though some exceptions exist, including those based on  intuitionistic logic (see \cite{Art02CSLI,ArtIem07JSL,Dashkov-JLC-2011,Steren-Bonelli-IMLA-2014,Studer-Marti-JAL-2016,Studer-Marti-AiML-2018,Kuznets-Marin-Strassburger-JAL-2021}),  on fuzzy logic (see \cite{Ghari-JIGPL-2016,Pischke-SL-2020,Pischke-IGPL-2021,Pischke-JIGPL-2023}), and on relevant logic (see \cite{Savic-Studer-2019,Standefe-IGPL-2019,Standefer2022,Standefer-JLLI-2023}).

To my knowledge, there has been no formalization of counterfactual (or subjunctive) conditionals in justification logics.  This paper aims to develop justification logics that is based on conditional logics, by incorporating expressions of the form  $\phi \limpliesCounterfactual \psi$, read as ``if it were that $\phi$, then it would be that $\psi$.''
The resulting logic, denoted by $\LPC$, is an extension of Artemov's logic of proofs $\Logic{LP}$ and Priest's conditional logic $\Logic{C}^+$ (see \cite{Priest2008}). We also explore   other conditional justification logics, including the logic $\LPCK$ that is an extension of the normal conditional logic $\Logic{CK}$ of Chellas \cite{Chellas1975}.

The semantics of conditional justification logics combines  selection-function semantics from conditional logics (see \cite{Chellas1975,Chellas-1980})  with  relational semantics from justification logics (see \cite{Sedlar-Podrouzek-LogKCA-2010}). Rather than using Chellas' selection-function, we employ   binary relations on states associated with each formula. (This equivalent approach is used, for instance, in \cite{Nute-Cross-HandbookPhilLogic-2001,Priest2008}.) We prove the completeness theorem of these logics via a canonical model construction using maximal consistent sets.

Counterfactual conditionals play an important role in many disciplines. For instance, in epistemology Nozick (1981, p. 179) defines knowledge as true belief satisfying certain conditions expressed by means of counterfactuals. As an application we formalize Nozick's counterfactual conditions, the sensitivity and adherence conditions, and investigate their relationship with  Aumman knowledge (i.e.  knowledge represented by $\Logic{S5}$-modalities with partition models). Additionally, we analyze  Gettier counterexamples to the justified true belief analysis and McGinn's counterexamples to Nozick's analysis of knowledge.

Relevant counterfactuals have also been examined in the literature (see \cite{Pereira-Aparicio-1989,Mares-Topoi-1994,Mares-Fuhrmann-1995,Mares-Book-2004,Priest2008}). We follow the presentation of Priest in \cite{Priest2008}[Section 10.7]. He adds the Routley star to the relational models of conditional logics and applies some modifications to the relational models to produce a logic with a relevant counterfactual. His logic is an extension of the relevant logic $\Logic{B}$ and conditional logic $\Logic{C}^+$. Following Priest’s presentation,  we present a conditional justification logic, called $\JRC$, that features a relevant counterfactual conditional. We show that the variable-sharing property holds for valid conditionals. We present a tableau system for $\JRC$ and show that it is complete with respect to Routly relational models. Finally, by formalizing the Gettier-type example of Chisholm, ``the sheep in the field'' from \cite{Chisholm-1966}, we show that neither the justified true belief nor Nozick's definition of knowledge is Aumann knowledge. 

\section{Language}
\label{sec:Syntax LPC}

The language of $\LPC$ is an extension of the language of $\LP$ with counterfactual conditionals. Let $\CTerms$ be a countable set of justification constants, $\VTerms$ a countable set of justification variables, and $\Prop$ a countable set of atomic propositions.

Justification terms and formulas are constructed by the following grammars:  
\begin{gather*}
	t \coloncolonequals \justConsThree \mid \justVarOne \mid  t \cdot t \mid  t + t \mid \ \tinspect t \, , \\
	\phi \coloncolonequals p \mid \neg \phi \mid \phi \wedge \phi \mid \phi \limplies \phi \mid \phi \limpliesCounterfactual \phi \mid \jbox{t} \phi \, .
\end{gather*}
where $\justConsThree \in \CTerms$, $\justVarOne \in \VTerms$, and $p \in \Prop$. The set of justification terms and formulas are denoted by  $\Terms$ and $\Formulae$, respectively. Formulas of the form $\jbox{t} \phi$ are called \textit{justification assertions} and are read ``the agent has justified true belief that $\phi$ justified by reason $t$.'' The conditional $\phi \limpliesCounterfactual \psi$ is read ``if it were that $\phi$, then it would be that $\psi$.''

We use the following usual abbreviations:
\begin{align*}
	\bot &\colonequals \phi \wedge \neg \phi, \\
\top &\colonequals \neg \bot, \\
\phi \lor \psi &\colonequals \neg (\neg \phi \wedge \neg \psi), \\
\phi \liff \psi &\colonequals (\phi \limplies \psi) \land (\psi \limplies \phi), \\
\phi \Leftrightarrow \psi &\colonequals (\phi \limpliesCounterfactual \psi) \land (\psi \limpliesCounterfactual \phi). 
\end{align*}

Associativity and precedence and the corresponding omission of brackets are handled as usual.

\section{Axiomatic system $\LPC$}
\label{sec:Axioms LPC}

We now begin with describing the axiom schemes and rules of our logic. The set of axiom schemes of $\LPC$ are:
\begin{enumerate}
	
	\item  All tautologies of classical propositional logic in the language of $\LPC$,
	
	\item  $(\phi \limpliesCounterfactual (\psi \limplies \chi)) \limplies ((\phi \limpliesCounterfactual \psi) \limplies (\phi \limpliesCounterfactual \chi))$,
	
	\item $\phi \limpliesCounterfactual \phi$.
	
	\item $(\phi \limpliesCounterfactual \psi) \limplies (\phi \limplies \psi)$.
	
	\item $(\jbox{s} (\phi \limpliesCounterfactual \psi) \wedge \jbox{t} \phi) \limpliesCounterfactual \jbox{s \cdot t} \psi$.
	
	\item  $\jbox{s} \phi \limpliesCounterfactual  \jbox{s \ttsum t} \phi$.
	
	\item  $\jbox{t} \phi \limpliesCounterfactual  \jbox{s \ttsum t} \phi$.
	
	\item $\jbox{t} \phi \limpliesCounterfactual \phi$. 
	
	\item $\jbox{t} \phi \limpliesCounterfactual \jbox{\tinspect t} \jbox{t} \phi$.
\end{enumerate}

The rules of inference are
\[
\lrule{ \phi \qquad  \phi \limplies \psi}{ \psi}\,\mprule 
\qquad
\lrule{  \psi}{ \phi \limpliesCounterfactual \psi}\, \RCNrule
\]

Axioms 2, 3, and 4 are denoted respectively by CK, ID, and MP in \cite{Chellas1975}. Axioms 5--9 have not been appeared in the literature of justification logic, though their classical version, where the conditional $\limpliesCounterfactual$ is replaced by $\limplies$, are standard (see e.g. \cite{Artemov-Fitting-Studer-SEP}).
 
A \textit{constant specification} $\CS$ for $\LPC$ is a set of formulas of the form
$\jbox{\justConsThree}  \phi$, 
where $\phi$ is an axiom instance of $\LPC$. Given a constant specification $\CS$ for $\LPC$, the justification logic $\LPC_\CS$ is defined as follows. The set of axioms of $\LPC_\CS$ is the above set of axioms 1--8 together with formulas from $\CS$, and the rules of inference are $\mprule$ and $\RCKrule$. Whenever we write $\LPC_\CS$, it means that  $\CS$ is a constant specification for $\LPC$. By $\vdash_{\LPC_\CS}$ we denote derivability  in $\LPC_\CS$.
We define a derivation from a set of premises $T$ as follows:
\[
T \vdash_{\LPC_\CS} \phi \text{ iff there exists $\psi_1,\ldots,\psi_n \in T$ such that 
	$\vdash_{\LPC_\CS} (\psi_1 \land \cdots \land \psi_n) \limplies \phi$.}
\]

From the above definition, it is obvious that the Deduction Theorem holds, i.e. if $T, \phi \vdash_{\LPC_\CS} \psi$, then $T \vdash_{\LPC_\CS} \phi \limplies \psi$.

We now show that the following rule is derivable in $\LPC$:
\[
\lrule{ (\psi_1 \wedge \cdots \wedge \psi_n) \limplies \psi}{ (\phi \limpliesCounterfactual \psi_1) \wedge \cdots \wedge (\phi \limpliesCounterfactual \psi_n) \limplies (\phi \limpliesCounterfactual \psi)}\, \RCKrule \quad (n \geq 0)
\]
Note that the rule $\RCNrule$ is a special case of $\RCKrule$ when $n = 0$. This rule is introduced by Chellas \cite{Chellas1975} in the axiom system of (normal) conditional logics. We will use this rule in the completeness proof of $\LPC_\CS$. We first need a lemma.

\begin{lemma}\label{lem: proof of axiom CC in LPC}
	Given any constant specification $\CS$ for $\LPC$, the following formula is provable in $\LPC_\CS$:
	\[
	((\phi \limpliesCounterfactual \psi_1) \wedge (\phi \limpliesCounterfactual \psi_2)) \limplies (\phi \limpliesCounterfactual (\psi_1 \wedge \psi_2)).
	\]
\end{lemma}
\begin{proof}
	The proof is as follows:
	\begin{enumerate}
		\item $\psi_1 \limplies (\psi_2 \limplies \psi_1 \wedge \psi_2)$ \hfill tautology
		
		\item $\phi \limpliesCounterfactual (\psi_1 \limplies (\psi_2 \limplies \psi_1 \wedge \psi_2))$ \hfill from 2 by $\RCNrule$
		
		\item $(\phi \limpliesCounterfactual \psi_1) \limplies (\phi \limpliesCounterfactual (\psi_2 \limplies \psi_1 \wedge \psi_2))$ \hfill from 2 and axiom 2, by $\mprule$
		
		\item $(\phi \limpliesCounterfactual (\psi_2 \limplies \psi_1 \wedge \psi_2)) \limplies ((\phi \limpliesCounterfactual \psi_2) \limplies (\phi \limpliesCounterfactual (\psi_1 \wedge \psi_2)))$ \hfill axiom 2
		
		\item $(\phi \limpliesCounterfactual \psi_1) \limplies ((\phi \limpliesCounterfactual \psi_2) \limplies (\phi \limpliesCounterfactual (\psi_1 \wedge \psi_2)))$ \hfill from 3 and 4 by propositional reasoning
		
		\item $(\phi \limpliesCounterfactual \psi_1) \wedge (\phi \limpliesCounterfactual \psi_2) \limplies (\phi \limpliesCounterfactual (\psi_1 \wedge \psi_2))$ \hfill from 5 by propositional reasoning \qed
	\end{enumerate}

\end{proof}

\begin{theorem}
	The rule $\RCKrule$ is derivable in $\LPC_\CS$.
\end{theorem}
\begin{proof}
	The proof is as follows:
	\begin{enumerate}
		\item $(\psi_1 \wedge \cdots \wedge \psi_n) \limplies \psi$ \hfill hypothesis
		
		\item $\phi \limpliesCounterfactual ((\psi_1 \wedge \cdots \wedge \psi_n) \limplies \psi)$ \hfill from 1 by $\RCNrule$
		
		\item $(\phi \limpliesCounterfactual (\psi_1 \wedge \cdots \wedge \psi_n)) \limplies (\phi \limpliesCounterfactual \psi))$ \hfill from 2 and axiom 2, by $\mprule$
		
		\item $(\phi \limpliesCounterfactual \psi_1) \wedge \cdots \wedge (\phi \limpliesCounterfactual \psi_n) \limplies (\phi \limpliesCounterfactual (\psi_1 \wedge \cdots \wedge \psi_n))$ \hfill by Lemma \ref{lem: proof of axiom CC in LPC}
		
		\item $(\phi \limpliesCounterfactual \psi_1) \wedge \cdots \wedge (\phi \limpliesCounterfactual \psi_n) \limplies (\phi \limpliesCounterfactual \psi)$ \hfill from 3 and 4 by propositional reasoning \qed
	\end{enumerate}
\end{proof}
\section{Relational semantics for $\LPC$}

The semantics of $\LPC$ integrates the selection-function semantics from conditional logics (see \cite{Chellas1975,Chellas-1980})  with the relational semantics from justification logics (see \cite{Sedlar-Podrouzek-LogKCA-2010}). Rather than using Chellas' selection-function, we utilize binary relations on states that are associated with each formula. (This equivalent approach is used, for instance, in \cite{Nute-Cross-HandbookPhilLogic-2001,Priest2008}.)

\begin{definition}\label{def: Relational models LPC}
	A relational model is a tuple $\model= (\worlds, \worldsNormal, \accrelFormula, \accrelTrem,  \valuation)$ where
	
	\begin{enumerate}
		\item $W$ and $W_N$ are non-empty sets of states and normal states, respectively, such that \linebreak $W_N \subseteq W$.

		\item $\accrelFormula$ is  a function that assigns to each formula $\phi \in \Formulae$ a binary relation $\relFormula$ on $\worldsNormal$. Given $w \in \worldsNormal$, let $\relFormula(w) := \{ v \in \worldsNormal \mid w \relFormula v \}$. 
		
		\item $\accrelTrem$ is  a function that assigns to each term $t \in \Terms$ a  binary relation $\relTrem$ on $\worlds$. Given $w \in \worlds$, let $\relTrem(w) := \{ v \in \worlds \mid w \relTrem v \}$. 
		
		\item $\valuation$ is a valuation function such that $\valuation: \worldsNormal \to \powerset(\Prop)$ and $\valuation: \worldsNonNormal \to \powerset(\Formulae)$. 	
	\end{enumerate}
\end{definition}

Here, $w \relFormula v$ means that ``$\phi$ is true at $v$, which is, ceteris paribus, the
same as $w$.'' In other words, $\relFormula(w)$ contains the most similar states to $w$ where $\phi$ is true. $w \relTrem v$ means that  ``at the state $w$, $v$ is doxastically possible for the agent according to reason $t$,'' and $\relTrem(w)$ can be viewed as doxastic alternatives of $w$ according to reason $t$.

\begin{definition}\label{def: Truth LPC}
	Let $\model= (\worlds, \worldsNormal, \accrelFormula, \accrelTrem,  \valuation)$ be a relational model. 
	Truth at non-normal states $w \in \worldsNonNormal$ are defined as follows:
	\[
	(\model, w) \entails \phi \quad \mbox{ iff } \quad \phi \in \valuation(w).
	\]
	Truth at normal states $w \in \worldsNormal$ is inductively defined as follows:
	\begin{align*}
		(\model, w) &\entails p \text{ iff } p \in \valuation(w) \, ,\\
		(\model, w) &\entails \neg \phi \text{ iff } (\model, w) \not\entails \phi \, ,\\
		(\model, w) &\entails \phi \wedge \psi \text{ iff } (\model, w) \entails \phi \text{ and } (\model, w) \entails \psi \, ,\\
		(\model, w) &\entails \phi \limplies \psi \text{ iff } (\model, w) \entails \phi \text{ implies } (\model, w) \entails \psi \, ,\\
		(\model, w) &\entails \phi \limpliesCounterfactual \psi \text{ iff } \relFormula(w) \subseteq \truthsetModel{\psi} \, ,\\
		(\model, w) &\entails \jbox{t} \phi \text{ iff } \relTrem(w) \subseteq \truthsetModel{\phi} \, ,
	\end{align*}
	where $\truthsetModel{\phi} := \{ w \in \worlds \mid (\model, w) \entails \phi \}$ is the \textit{truth set} of $\phi$ in $\model$.
\end{definition}	
The truth conditions of counterfactual conditionals and justification assertions can be rewritten as follows:
\begin{align*}
(\model, w) &\entails \phi \limpliesCounterfactual \psi \text{ iff }   \text{ for all normal states }  v \in \worldsNormal\  (\text{if }   w \relFormula  v, \text{ then }  (\model, v) \entails \phi)\, , \\
(\model, w) &\entails \jbox{t} \phi \text{ iff }   \text{ for all states }  v \in \worlds\ (\text{if }   w \relTrem  v, \text{ then }  (\model, v) \entails \phi) \,.			
\end{align*}		
The above conditions show that $\phi \limpliesCounterfactual \psi$ and $\jbox{t} \phi$ can be viewed as necessity modalities. One may consider $\phi \limpliesCounterfactual \psi$ as the relative necessity $\jbox{\phi} \psi$, expressing that $\psi$ is necessary relative to $\phi$ (cf. \cite{Chellas1975}).

\begin{definition}\label{def: relational models conditions LPC}
	A relational model $\model= (\worlds, \worldsNormal, \accrelFormula, \accrelTrem,  \valuation)$ is called an $\LPC_\CS$-model if it satisfies the following conditions. For all $w \in \worldsNormal$:
	\begin{enumerate}
		
		\item For all $\phi \in \Formulae$, $\relFormula(w) \subseteq \truthsetModel{\phi}$.
		
		\item For all $\phi \in \Formulae$, if $w \in \truthsetModel{\phi}$, then $w \in \relFormula(w)$.
		
		\item For  all $\jbox{\justConsThree}  \phi \in \CS$,	$\accrel_{\justConsThree}(w) \subseteq \truthsetModel{\phi}$.
		
		\item For all $s,t \in \Terms$, $\accrel_{s \ttsum t}(w) \subseteq \accrel_{s}(w) \cap \accrel_{t}(w)$.

		\item For all $s,t \in \Terms$,	 $\accrel_{s \cdot t} (w)  \subseteq (\accrel_{s}  \times \accrel_{t}) (w)$, where
		\begin{multline*}
			\accrel_{s}  \times \accrel_{t}  := \{ (w,v) \in \worlds\times \worlds \mid \text{for all $\phi, \psi \in \Formulae$: } \\
			\text{ if  } w \in \truthsetModel{\jbox{s} (\phi \limpliesCounterfactual \psi) \wedge \jbox{t} \phi}, \text{ then } v \in \truthsetModel{\psi} \};
		\end{multline*}
		

		\item For all $t \in \Terms$, $w \accrel_{t} w$.
		
		\item For all $t \in \Terms$ and all $v, u \in \worlds$, if $w \accrel_{!t} v$ and $v \accrel_{t} u$, then $w \accrel_{t} u$.
	\end{enumerate}
\end{definition}

\begin{definition}\label{def: validity LPC}
	Let $\CS$ be a constant specification for $\LPC$.
	
	\begin{enumerate}
		\item Given an $\LPC_\CS$-model $\model= (\worlds, \worldsNormal, \accrelFormula, \accrelTrem,  \valuation)$, we write $\model \entails \phi$ if
		for all $w \in \worldsNormal$, we have 
		$(\model, w) \entails \phi$.
		
		\item A formula $\phi$ is $\LPC_\CS$-valid, denoted by $\entails_{\LPC_\CS} \phi$, if $\model \entails \phi$ for all $\LPC_\CS$-models $\model$.
		
		\item Given a set of formulas $T$ and a formula $\phi$ of $\LPC_\CS$, the (local) consequence relation is defined as follows:\\ $T \entails_{\LPC_\CS} \phi$ iff for all 
		$\LPC_\CS$-models $\model= (\worlds, \worldsNormal, \accrelFormula, \accrelTrem,  \valuation)$, for all $w \in \worldsNormal$, if $(\model, w) \entails \psi$ for all $\psi \in T$, then $(\model, w) \entails \phi$.
	\end{enumerate}
	
\end{definition}

%

\subsection*{Counterpossibles}

Counterpossibles, or counterfactuals with impossible antecedents, have sparked significant debate regarding their truth conditions. Lewis--Stalnaker and also Chellas semantics of conditional logic have it that all counterpossibles  are vacuously true
 (for a discussion on counterpossibles see e.g. \cite{Brogaard-Salerno-2013,Williamson-Topoi-2018,Berto-et-al-2018,McLoone-et-al-2023}). We show that the same holds in $\LPC$.

\begin{lemma}
	Given any constant specification $\CS$ for $\LPC$, for any $\phi \in \Formulae$ we have
	\begin{enumerate}
		\item $ \entails_{\LPC_\CS} \bot \limpliesCounterfactual \phi$.
		
		\item $ \entails_{\LPC_\CS} \phi \limpliesCounterfactual \top$.
	\end{enumerate}
\end{lemma}
\begin{proof}
	Suppose that $\phi \in \Formulae$ is an arbitrary formula, $\model= (\worlds, \worldsNormal, \accrelFormula, \accrelTrem,  \valuation)$ is an arbitrary $\LPC_\CS$-model, and  $w \in \worldsNormal$ is an arbitrary normal state. Then, since $\accrel_{\bot} (w) \subseteq \worldsNormal$ and $\truthsetModel{\bot} \subseteq \worldsNonNormal$, from condition 1 of Definition \ref{def: canonical relational model LPC} (i.e. $\accrel_{\bot} (w) \subseteq \truthsetModel{\bot}$), it follows that $\accrel_{\bot} (w) = \emptyset$. Hence, $\accrel_{\bot} (w) \subseteq \truthsetModel{\phi}$. Therefore, $(\model, w) \entails \bot \limpliesCounterfactual \phi$. Furthermore, since $\accrel_{\phi} (w) \subseteq \worldsNormal \subseteq \truthsetModel{\top}$, we get $ (\model, w) \entails \phi \limpliesCounterfactual \top$. Therefore, $ \entails_{\LPC_\CS} \bot \limpliesCounterfactual \phi$ and $ \entails_{\LPC_\CS} \phi \limpliesCounterfactual \top$.\qed
\end{proof}

In section \ref{sec: Tableaux for JRC}, we introduce a conditional justification logic in which counterpossibles are not valid (see Lemma \ref{lem: counterpossible in JCR}).
\subsection*{Hyperintentionality}

In normal modal logics, the following regularity rule is admissible:
\[
\lrule{ \phi \liff \psi}{ \Box \phi \liff \Box \psi}\ (\sf{RE})
\]
However, in justification logic from the fact that $\phi$ and $\psi$ are equivalent and that $\phi$ is justified by reason $t$, it is not necessarily concluded that $\psi$ is justified by the same reason $t$. In other words, justifications are hyperintentional (see  \cite{Art-Fit-Book-2019} and \cite{Faroldi-Protopopescu-IGPL-2019} for a more detailed exposition). We show that justification versions of the rule $(\sf{RE})$ is not admissible in the justification logic $\LPC$.

\begin{lemma}
	The regularity rule
	\[
	\lrule{\phi \liff \psi}{\jbox{t} \phi \liff \jbox{t} \psi}\ (\sf{JRE}_\equiv)
	\]
	is not validity preserving in $\LPC_\CS$, i.e. there exists $\phi, \psi \in \Formulae$ and $t \in \Terms$ such that \linebreak $\entails_{\LPC_\CS} \phi \liff \psi$ but $\not \entails_{\LPC_\CS} \jbox{t} \phi \liff \jbox{t} \psi$.
\end{lemma}
\begin{proof}
	Let $\CS$ be an arbitrary constant specification for $\LPC$. Consider the relational model $\model = (\worlds, \worldsNormal, \accrelFormula, \accrelTrem,  \valuation)$ as follows:
	\begin{itemize}
		\item $\worlds = \{ w, v\}$,
		
		\item $\worldsNormal = \{ w \}$,
		
		\item $\accrel_{\phi}(w) = \truthsetModel{\phi} \cap \worldsNormal$, for all $\phi \in \Formulae$,
		
		\item $\accrel_{\justVarOne}(w) = \{ w, v\}$, for all $\justVarOne \in \VTerms$, \\
				 $\accrel_{\justConsThree}(w) = \{ w \}$  for all $\justConsThree \in \CTerms$,\\
				$\accrel_{t}(w) = \{ w\}$, for all $t \in \Terms \setminus (\VTerms \cup \CTerms)$, \\
		$\accrel_{t}(v) = \emptyset$, for all $t \in \Terms$,
		
		\item $\valuation(w) = \{ p \}$, $\valuation(v) = \{ p \wedge p \}$.
	\end{itemize}
	It is easy to show that $\model$ is an $\LPC_\CS$-model. 
	
	Let $\justVarOne \in \VTerms$. We observe that $(\model,w) \not \entails \jbox{\justVarOne} p$ and $(\model,w)  \entails \jbox{\justVarOne} (p \wedge p)$, because \linebreak $\accrel_{\justVarOne}(w) = \{ w, v\} \not \subseteq \truthsetModel{p} = \{ w \}$ and $\accrel_{\justVarOne}(w) = \{ w, v\}  \subseteq\truthsetModel{p \wedge p} = \{w, v\}$. Hence,\linebreak $(\model,w) \not \entails \jbox{\justVarOne} p \liff \jbox{\justVarOne} (p \wedge p)$. Thus, $\not \entails_{\LPC_\CS} \jbox{\justVarOne} p \liff \jbox{\justVarOne} (p \wedge p)$, but $\entails_{\LPC_\CS} p \liff p \wedge p$.  \qed
\end{proof}

We will now demonstrate that the regularity rule remains inadmissible, even when the material biconditional $\equiv$ in its formulation is replaced with the counterfactual biconditional $\Leftrightarrow$.

\begin{lemma}
	The regularity rule
	\[
	\lrule{\phi \Leftrightarrow \psi}{\jbox{t} \phi \Leftrightarrow \jbox{t} \psi}\ (\sf{JRE}_\Leftrightarrow)
	\]
	is not validity preserving in $\LPC_\CS$, i.e. there exists $\phi, \psi \in \Formulae$ and $t \in \Terms$ such that \linebreak $\entails_{\LPC_\CS} \phi \Leftrightarrow \psi$ but $\not \entails_{\LPC_\CS} \jbox{t} \phi \Leftrightarrow \jbox{t} \psi$.
\end{lemma}

\begin{proof}
	Let $\CS$ be an arbitrary constant specification for $\LPC$. Consider the relational model $\model = (\worlds, \worldsNormal, \accrelFormula, \accrelTrem,  \valuation)$ as follows:
	\begin{itemize}
		\item $\worlds =  \{ w, u, v\}$,
		
		\item $\worldsNormal = \{ w, u \}$,
		
		\item $\accrel_{\phi}(w) = \accrel_{\phi}(u) =\truthsetModel{\phi} \cap \worldsNormal$, for all $\phi \in \Formulae$,
		
		\item $\accrel_{\justVarOne}(w) = \{ w\}$, for all $\justVarOne \in \VTerms$, \\
		$\accrel_{\justVarOne}(u) = \{ u, v\}$, for all $\justVarOne \in \VTerms$, \\
		$\accrel_{\justConsThree}(w) = \accrel_{\justConsThree}(u) = \{ w, u \}$  for all $\justConsThree \in \CTerms$,\\
		$\accrel_{t}(w) = \{ w \}$, $\accrel_{t}(u) = \{  u \}$ for all $t \in \Terms \setminus (\VTerms \cup \CTerms)$, \\
		$\accrel_{t}(v) = \emptyset$, for all $t \in \Terms$,
		
		\item $\valuation(w) = \valuation(u) = \{ p \}$, $\valuation(v) = \{  p \wedge p \limpliesCounterfactual p \}$.
	\end{itemize}
	It is easy to show that $\model$ is an $\LPC_\CS$-model. 
	
	Let $\justVarOne \in \VTerms$. We observe that $(\model,w) \not \entails  \jbox{x} (p  \wedge p \limpliesCounterfactual p) \limpliesCounterfactual \jbox{\justVarOne} ( p \vee\neg p)$, because 
	$$\accrel_{\jbox{\justVarOne} (p \wedge p \limpliesCounterfactual p)} (w) = \truthsetModel{\jbox{\justVarOne} (p \wedge p \limpliesCounterfactual p)}  \cap \worldsNormal = \{ w,u \} \not \subseteq \truthsetModel{\jbox{\justVarOne} ( p \vee\neg p)} = \{w\}.$$
	Thus, $\not \entails_{\LPC_\CS} \jbox{\justVarOne} (p \wedge p \limpliesCounterfactual p) \Leftrightarrow \jbox{\justVarOne} (p \vee \neg p)$, but $\entails_{\LPC_\CS}  (p  \wedge p \limpliesCounterfactual p) \Leftrightarrow (p \vee \neg p)$.  \qed
\end{proof}

\section{Completeness of $\LPC$}
\label{sec: Completeness of LPC}

For the completeness proof, we employ the canonical model construction. Given a logic $L$ and a set of formulas $T$ of $L$, we say that $T$ is $L$-consistent if $T \not \vdash_{L} \bot$. We say that $T$ is maximal $L$-consistent if it is $L$-consistent and it has no proper $L$-consistent extension. Any $L$-consistent set can be extended to a maximal $L$-consistent set by systematically adding each formula that preserves consistency  (indeed, a Lindenbaum-like construction can be used to show this). The set of all maximal $L$-consistent sets of formulas is denoted by $\MCS_{L}$.

Note that from axiom 3 it follows that if $\Gamma \in \MCS_{\LPC_\CS}$ and $\phi \limpliesCounterfactual \psi, \phi \in \Gamma$, then $\psi \in \Gamma$ (in fact, since $\Gamma$ is a maximal $\LPC_\CS$-consistent set, $(\phi \limpliesCounterfactual \psi) \limplies (\phi \limplies \psi)\in \Gamma$, and thus $\phi \limplies \psi \in \Gamma$). This fact will be used several times in this section.

The following auxiliary lemma is used in the completeness proof.
\begin{lemma}\label{lem: Fact in completeness of LPC}
	Let $T$ be a set of formulas and $\phi \in \Formulae$. If $\phi \in \cap \lset{\Delta  \in \MCS_{\LPC_\CS} \mid T \subseteq \Delta}$, then $T \vdash_{\LPC_\CS} \phi$.
\end{lemma}
\begin{proof}
	Assume that $\phi \in \cap \lset{\Delta  \in \MCS_{\LPC_\CS} \mid T \subseteq \Delta}$, and suppose towards a contradiction that $T \not \vdash_{\LPC_\CS} \phi$. Then, $T \cup \lset{ \neg \phi }$ is $\LPC_\CS$-consistent and there is $\Delta  \in \MCS_{\LPC_\CS}$ such that $T \cup \lset{ \neg \phi } \subseteq \Delta$.  On the other hand, from the assumption we have that $\phi \in \Delta$, which is a contradiction. \qed
\end{proof}

\begin{definition}\label{def: canonical relational model LPC}
	The canonical relational model 
	\[
	\model = (\worlds, \worldsNormal, \accrelFormula, \accrelTrem,  \valuation)
	\] 
	for $\LPC_\CS$ is defined as follows:
	\begin{enumerate}
		\setlength\itemsep{0.1cm}	
		\item 
		$\worlds := \powerset(\Formulae) = \lset{ \Gamma \mid \Gamma \subseteq \Formulae}$;
		
		\item 
		$\worldsNormal \colonequals \MCS_{\LPC_\CS}$;
		
		\item 
		$\Gamma \relFormula \Delta$ iff $\Gamma / \phi \subseteq \Delta$, for $\Gamma, \Delta \in \worldsNormal$, 
		where $\Gamma /\phi = \{ \psi \in \Formulae \mid \phi \limpliesCounterfactual \psi \in \Gamma \}$;
		
		\item
		$\Gamma \relTrem \Delta$ iff $\Gamma / t \subseteq \Delta$, for $\Gamma, \Delta \in \worlds$, 
		where $\Gamma /t = \{ \psi \in \Formulae \mid \jbox{t} \psi \in \Gamma \}$;
		
		\item 
		$\valuation(\Gamma) \colonequals \Gamma$, for $\Gamma \in \worldsNonNormal$;\\
		$\valuation(\Gamma) \colonequals \Prop \cap \Gamma$, for $\Gamma \in \worldsNormal$.
	\end{enumerate}
	
\end{definition}

\begin{lemma}[Truth Lemma]\label{lem: Truth Lemma LPC}
	Let 
	\(
	\model= (\worlds, \worldsNormal, \accrelFormula, \accrelTrem,  \valuation)
	\) 
	be the canonical relational model for $\LPC_\CS$. Then for any $\phi \in \Formulae$ and any $\Gamma \in \worlds$
	\[
	(\model, \Gamma) \entails \phi \quad \mbox{ iff }  \quad \phi \in \Gamma.
	\] 
\end{lemma}
\begin{proof}
	By induction on $\phi$. If  $\Gamma \in \worldsNonNormal$ we have
	\[
	(\model, \Gamma) \entails \phi \quad \mbox{ iff } \quad  \phi \in \valuation(\Gamma) = \Gamma.
	\]	
	Suppose $\Gamma \in \worldsNormal$. For $\phi = p \in \Prop$ we have
	\[
	(\model, \Gamma) \entails p \quad \mbox{ iff } \quad p \in \valuation(\Gamma) \quad \mbox{ iff } \quad p \in \Gamma.
	\]	
	The cases in which the main connective of $\phi$ is either $\neg, \wedge$, or $\limplies$ are standard. We only check the cases $\phi = \phi_1 \limpliesCounterfactual \phi_2$ and $\phi = \jbox{t} \psi$.
	
	If $(\model, \Gamma) \entails \phi_1 \limpliesCounterfactual \phi_2$, then $\relModel{\phi_1} (\Gamma) \subseteq \truthsetModel{\phi_2}$. Let $T := \Gamma / \phi_1 =  \{ \psi \mid \phi_1 \limpliesCounterfactual \psi \in \Gamma \}$. Given an arbitrary $\Delta \in \worldsNormal$ such that $T \subseteq \Delta$, we have that $(\model, \Delta) \entails \phi_2$, and thus by the induction hypothesis we get $\phi_2 \in \Delta$. Thus, 
	\[
	\phi_2 \in \cap \lset{\Delta  \in \MCS_{\LPC_\CS} \mid T \subseteq \Delta}.
	\] 
	By Lemma  \ref{lem: Fact in completeness of LPC}, we have that $T \vdash_{\LPC_\CS} \phi_2$. Thus, $\vdash_{\LPC_\CS} \psi_1 \wedge \ldots \wedge \psi_n \limplies \phi_2$, for some $\psi_1, \ldots, \psi_n \in T$. Therefore, $\psi_1 \wedge \ldots \wedge \psi_n \limplies \phi_2 \in \Gamma$ and also $\phi_1 \limpliesCounterfactual \psi_1, \ldots, \phi_1 \limpliesCounterfactual \psi_n \in \Gamma$. Using the derived rule $\RCKrule$, we get $\phi_1 \limpliesCounterfactual \phi_2 \in \Gamma$, as required. 
	
	Now suppose $\phi_1 \limpliesCounterfactual \phi_2 \in \Gamma$. Then, if $\Gamma \relModel{\phi_1} \Delta$, for an arbitrary $\Delta \in \worldsNormal$, then $\Gamma / \phi_1 \subseteq \Delta$, and thus $\phi_2 \in \Delta$.  By the induction hypothesis, $(\model, \Delta) \entails \phi_2$. Therefore,  $(\model, \Gamma) \entails \phi_1 \limpliesCounterfactual \phi_2$.
	
	If $(\model, \Gamma) \entails \jbox{t} \psi$, then by letting $\Delta := \Gamma /t$ we get $\Gamma \relModel{t} \Delta$, and thus $(\model, \Delta) \entails \psi$. By the induction hypothesis, $\psi \in \Delta$, and hence $\jbox{t} \psi \in \Gamma$.
	
	Suppose $\jbox{t} \psi \in \Gamma$. Then, if $\Gamma \relModel{t} \Delta$, for an arbitrary $\Delta \in \worlds$, then $\Gamma /t \subseteq \Delta$, and thus $\psi \in \Delta$.  By the induction hypothesis, $(\model, \Delta) \entails \psi$. Therefore,  $(\model, \Gamma) \entails \jbox{t} \psi$.
	\qed
\end{proof}

\begin{lemma}\label{lem: canonical model is a relational model LPC}
	The canonical relational model
	\(
	\model= (\worlds, \worldsNormal, \accrelFormula, \accrelTrem,  \valuation)
	\) 
	 is an $\LPC_\CS$-model. 
\end{lemma}
\begin{proof}
	It suffices to show that the canonical model $\model= (\worlds, \worldsNormal, \accrelFormula, \accrelTrem,  \valuation)$ satisfies conditions 1--7 of Definition \ref{def: relational models conditions LPC}. 
	Let $\Gamma \in \worldsNormal =  \MCS_{\LPC_\CS}$.  
	\begin{enumerate}
		\item Given $\phi \in \Formulae$, suppose that $\Delta$ is an arbitrary member of $\relFormula (\Gamma)$. Then, $\Gamma / \phi \subseteq \Delta$. Since $\Gamma$ is a maximal $\LPC_\CS$-consistent set, $\phi \limpliesCounterfactual \phi \in \Gamma$. Thus, $\phi \in \Delta$, and hence, by the Truth Lemma, $\Delta \in \truthsetModel{\phi}$. Therefore, $\relFormula (\Gamma) \subseteq \truthsetModel{\phi}$.
		
		\item Given $\phi \in \Formulae$, suppose that $\Gamma \in \truthsetModel{\phi}$. Hence, by the Truth Lemma, $\phi \in \Gamma$. We want to show that $\Gamma / \phi \subseteq \Gamma$. Suppose that $\psi \in \Gamma / \phi$, and hence $\phi \limpliesCounterfactual \psi \in \Gamma$.  Therefore, $\psi \in \Gamma$.
		
		\item Suppose that  $\jbox{\justConsThree} \phi \in \CS$, and $\Delta$ is an arbitrary member of $\relModel{\justConsThree} (\Gamma)$. Then, $\Gamma / \justConsThree \subseteq \Delta$. Since $\Gamma$ is a maximal $\LPC_\CS$-consistent set, $\jbox{\justConsThree} \phi \in \Gamma$. Thus, $\phi \in \Delta$, and hence, by the Truth Lemma, $\Delta \in \truthsetModel{\phi}$. Therefore, $\relModel{\justConsThree} (\Gamma) \subseteq \truthsetModel{\phi}$.
		
		\item Given $s,t \in \Terms$, suppose that $\Delta$ is an arbitrary member of $\relModel{s+t} (\Gamma)$. We want to show that $\Gamma / s, \Gamma / t \subseteq \Delta$. Let $\phi \in \Gamma / s$. Then, $\jbox{s} \phi \in \Gamma$. Since $\Gamma$ is a maximal $\LPC_\CS$-consistent set, we have that $\jbox{s} \phi \limpliesCounterfactual \jbox{s+t} \phi \in \Gamma$, and thus $\jbox{s+t} \phi \in \Gamma$. Since $\Gamma / (s+ t) \subseteq \Delta$, we get $\phi \in \Delta$. Therefore, $\Gamma / s \subseteq \Delta$. Similarly, it can be shown that $\Gamma /  t \subseteq \Delta$.
		
		\item Given $s,t \in \Terms$, suppose that $\Delta$ is an arbitrary member of $\relModel{s \cdot t} (\Gamma)$. We want to show that $\Delta \in (\relModel{t} \times \relModel{t}) (\Gamma)$. Let $\Gamma \in \truthsetModel{\jbox{s} (\phi \limpliesCounterfactual \psi) \wedge \jbox{t} \phi}$, for some arbitrary $\phi, \psi \in \Formulae$. We have to show that $\Delta \in \truthsetModel{\psi}$. By the Truth Lemma, $\jbox{s} (\phi \limpliesCounterfactual \psi) \wedge \jbox{t} \phi \in \Gamma$. Since $\Gamma$ is a maximal $\LPC_\CS$-consistent set, we have that $\jbox{s} (\phi \limpliesCounterfactual \psi) \wedge \jbox{t} \phi \limpliesCounterfactual \jbox{s \cdot t} \psi \in \Gamma$, and hence $\jbox{s \cdot t} \psi \in \Gamma$. Since $\Gamma / (s \cdot t) \subseteq \Delta$, we get $\psi \in \Delta$. Therefore, by the Truth Lemma, $\Delta \in \truthsetModel{\psi}$.
		
		\item Given $t \in \Terms$, we want to show that $\Gamma \relTrem \Gamma$. Suppose that $\phi \in \Gamma /t$, and hence $\jbox{t} \phi \in \Gamma$. Since $\Gamma$ is a maximal $\LPC_\CS$-consistent set, $\jbox{t} \phi \limpliesCounterfactual \phi \in \Gamma$, and hence $\phi \in \Gamma$. Therefore, $\Gamma / t \subseteq \Gamma$. 
		
		\item Given $t \in \Terms$, suppose that $\Gamma \relModel{!t} \Delta$ and $\Delta \relModel{t} \Sigma$, for $\Delta , \Sigma \in \worlds$. We want to show that $\Gamma \relModel{t} \Sigma$ or $\Gamma / t \subseteq \Sigma$. Suppose that $\phi \in \Gamma / t$, and hence $\jbox{t} \phi \in \Gamma$. Since $\Gamma$ is a maximal $\LPC_\CS$-consistent set, we have that $\jbox{t} \phi \limpliesCounterfactual \jbox{!t} \jbox{t} \phi \in \Gamma$, and thus $\jbox{!t} \jbox{t} \phi \in \Gamma$. Since $\Gamma / !t \subseteq \Delta$, we get $\jbox{t} \phi \in \Delta$. Since $\Delta / t \subseteq \Sigma$, we get $\phi \in \Sigma$,  as required. \qed
	\end{enumerate}
\end{proof}

Now completeness follows from Lemmas \ref{lem: Truth Lemma LPC} and \ref{lem: canonical model is a relational model LPC}.

\begin{theorem}[Completeness]
	Let $\CS$ be a constant specification for $\LPC$. For each formula $\phi$ and set of formulas $T$,
	\[
	T \entails_{\LPC_\CS}  \phi  \quad\text{iff}\quad  T \vdash_{\LPC_\CS} \phi.
	\]
\end{theorem}
\begin{proof}
	Soundness is straightforward. For completeness, suppose that $T \not \vdash_{\LPC_\CS} \phi$. Thus,\linebreak $T \cup \{ \neg \phi \}$ is $\LPC_\CS$-consistent, and hence there is $\Gamma \in \MCS_{\LPC_\CS}$ such that $T \cup \{ \neg \phi \} \subseteq \Gamma$. In the canonical  model
	\(
	\model= (\worlds, \worldsNormal, \accrelFormula, \accrelTrem,  \valuation)
	\) 
	for $\LPC_\CS$, we have that $\Gamma \in \worldsNormal$. By the Truth Lemma~\ref{lem: Truth Lemma LPC}, we get  $\Gamma \in \truthsetModel{\neg \phi}$ and $\Gamma \in \truthsetModel{\psi}$ for all $\psi \in T$. Hence, $T \not \entails_{\LPC_\CS} \phi$. \qed
\end{proof}

\section{Other axioms and rules}
\label{sec: Other axioms and rules}

\subsection{The logic $\LPCint$}
\label{sec: The logic LPCint}

Justification logics are modal-like logics where justification terms replace modalities. The correspondence between modal and justification logics is established by the so-called ``realization theorem'' (cf. \cite{Art01BSL,Fit16APAL}), which specifies how occurrences of modalities in a modal formula are replaced by justification terms so that theorems of modal logic are transformed into theorems of justification logic. In order to show this correspondence, a result similar to the necessitation rule of normal modal logics (every theorem is necessitated) is needed in the context of justification logics. The ``internalization lemma'' plays such a role. It states that every theorem is justified. In this section, we extend the language and axioms of $\LPC$ to a logic called $\LPCint$ that enjoys the internalization property.

Justification terms and formulas of $\LPCint$ are constructed by the following mutual grammar:  
\begin{gather*}
	t \coloncolonequals \justConsThree \mid \justVarOne \mid  t \cdot t \mid  t + t \mid \ \tinspect t \mid \langle t, \phi \rangle  \, , \\
	\phi \coloncolonequals p \mid \neg \phi \mid \phi \wedge \phi \mid \phi \limplies \phi \mid \phi \limpliesCounterfactual \phi  \mid \jbox{t} \phi \, .
\end{gather*}
For example, the following are justification terms of $\LPCint$: 
$$\langle \tinspect \justConsThree, \neg p  \rangle, ~~\langle \justVarOne \cdot \justVarTwo, \jbox{\justConsThree + \justConsThree} (p \limpliesCounterfactual q) \rangle.$$

Let $\LPCint$ denote the logic obtained from $\LPC$ by adding the following axiom:
\begin{description}
	\item[$10.$] $\jbox{t}  \psi \limplies \jbox{\langle t, \phi \rangle} (\phi \limpliesCounterfactual \psi)$.
\end{description}
The definition of constant specifications for $\LPCint$, and the definition of $\LPCint_\CS$ and  derivability in $\LPCint_\CS$, denoted by $\vdash_{\LPCint_\CS}$, is similar to that of $\LPC_\CS$. 

Relational models of $\LPCint$ are defined similar to $\LPC$-models satisfying an additional condition corresponding to axiom 10. An $\LPCint_\CS$-model $\model = (\worlds, \worldsNormal, \accrelFormula, \accrelTrem,  \valuation)$ is an $\LPC_\CS$-model such that in addition to the conditions of Definition \ref{def: relational models conditions LPC} it satisfies the following extra condition:
\begin{description}
	\item[$8.$] For all $w \in \worldsNormal$, for all $v \in \worlds$, for all $\phi, \psi \in \Formulae$, and for all $t \in \Terms$: if $w \in \truthsetModel{\jbox{t} \psi}$ and $w \accrel_{\langle t, \phi \rangle} v$, then $v \in \truthsetModel{\phi \limpliesCounterfactual \psi}.$
\end{description}
Validity is defined similar to Definition \ref{def: validity LPC} and is denoted by $\models_{\LPCint_\CS}$. We now show the completeness of $\LPCint$.
\begin{theorem}[Completeness]
	Let $\CS$ be a constant specification for $\LPCint$. For each formula $\phi$ and set of formulas $T$,
	\[
	T \entails_{\LPCint_\CS}  \phi  \quad\text{iff}\quad  T \vdash_{\LPCint_\CS} \phi.
	\]
\end{theorem}
\begin{proof}
The proof is similar to that given in Section \ref{sec: Completeness of LPC}. The major difference is that we have to take maximal $\LPCint_\CS$-consistent sets into consideration. The proof of the Truth Lemma is similar to that of Lemma \ref{lem: Truth Lemma LPC}. Thus, it is enough to show that the canonical relational model $\model = (\worlds, \worldsNormal, \accrelFormula, \accrelTrem,  \valuation)$ for $\LPCint_\CS$ is an $\LPCint_\CS$-model. It suffices to show that $\model$ satisfies condition 8.

Suppose that for $\Gamma \in \worldsNormal$, $\Delta \in \worlds$,  $\phi, \psi \in \Formulae$, and $t \in \Terms$ we have that $\Gamma \in \truthsetModel{\jbox{t} \psi}$ and $\Gamma \accrel_{\langle t, \phi \rangle} \Delta$. Since $\Gamma$ is a maximal $\LPCint_\CS$-consistent set, $\jbox{t}  \psi \limplies \jbox{\langle t, \phi \rangle} (\phi \limpliesCounterfactual \psi) \in \Gamma$. By the Truth Lemma, from $\Gamma \in \truthsetModel{\jbox{t} \psi}$ it follows that $ \jbox{t} \psi \in \Gamma$, and therefore  $\jbox{\langle t, \phi \rangle} (\phi \limpliesCounterfactual \psi) \in \Gamma$. From this and $\Gamma \accrel_{\langle t, \phi \rangle} \Delta$ we obtain $\phi \limpliesCounterfactual \psi \in \Delta$. By the Truth Lemma, we get $\Delta \in \truthsetModel{\phi \limpliesCounterfactual \psi}$ as desired. \qed
\end{proof}

Before we show the internalization property for $\LPCint$, we need a definition.

A constant specification $\CS$ is called \textit{axiomatically appropriate} provided for every axiom instance $\phi$  of $\LPCint$ there exists a proof constant $\justConsThree$ such that $\jbox{\justConsThree} \phi \in \CS$.

\begin{lemma}[Internalization]
	Given an axiomatically appropriate constant specification $\CS$ for $\LPCint$, if $\vdash_{\LPCint_\CS} \chi$, then $\vdash_{\LPCint_\CS} \jbox{t} \chi$, for some term $t$.
\end{lemma}
\begin{proof}
	The proof is by induction on the derivation of
	$\chi$. We have two base cases:
	\begin{itemize}
		\item If $\chi$ is an axiom instance of $\LPCint$. Then, since $\CS$ is axiomatically appropriate, there is a proof constant $\justConsThree$ such that $\jbox{\justConsThree} \chi \in \CS$. Thus, put $t:= \justConsThree$.
		
		\item If $\chi = \jbox{\justConsThree} \phi \in \CS$, then using axiom 8 and $\mprule$ we get $\vdash_{\LPCint_\CS} \jbox{!\justConsThree} \jbox{\justConsThree} \phi$. Thus, put $t:= !\justConsThree$.
	\end{itemize}
	For the induction step, we have again two cases:
	\begin{itemize}
		\item Suppose $\chi$ is obtained by $\mprule$ from $\phi \limplies \chi$ and $\phi$. By the induction hypothesis, there are terms $r$ and $s$ such that $\vdash_{\LPCint_\CS} \jbox{r}(\phi \limplies \chi)$ and $\vdash_{\LPCint_\CS} \jbox{s} \phi$. Then, put $t: = r \cdot s$ and use  axioms 3 and 4 and $\mprule$ to obtain $\vdash_{\LPCint_\CS} \jbox{r \cdot s} \chi$.
	
		\item Suppose $\chi = \phi \limpliesCounterfactual \psi$ is obtained by $\RCKrule$ from  $\psi$. By the induction hypothesis, there is a term $s$ such that $\vdash_{\LPCint_\CS} \jbox{s} \psi$. Using axiom 10, we get $\vdash_{\LPCint_\CS} \jbox{\langle s, \phi \rangle} \chi$. Thus, put $t:= \langle s, \phi \rangle$.
	\end{itemize}
	\qed	
\end{proof}

\subsection{The logic $\LPCstandard$}

The classical justification logics usually includes an application axiom scheme  of the form $\jbox{s} (\phi \limplies \psi) \limplies (\jbox{t} \phi \limplies \jbox{s \cdot t} \psi)$. In the case that the base logic is classical, the application axiom is equivalent to $(\jbox{s} (\phi \limplies \psi) \wedge \jbox{t} \phi) \limplies \jbox{s \cdot t} \psi$. However, since the two formulas $\phi \limpliesCounterfactual (\psi \limpliesCounterfactual \chi)$ and $\phi \wedge \psi \limpliesCounterfactual \chi$ are not equivalent in $\LPC$, the axiom 4 of  $\LPC$, i.e. $(\jbox{s} (\phi \limpliesCounterfactual \psi) \wedge \jbox{t} \phi) \limpliesCounterfactual \jbox{s \cdot t} \psi$, is not equivalent to the following axiom
\begin{description}
	\item[$4'$.] $\jbox{s} (\phi \limpliesCounterfactual \psi) \limpliesCounterfactual (\jbox{t} \phi \limpliesCounterfactual \jbox{s \cdot t} \psi)$.
\end{description}
%
In this section, we show that it is quite possible to axiomatize conditional justification logics using axiom $4'$ instead of axiom $4$.

Let $\LPCstandard$ denote the logic obtained from $\LPC$ by replacing axiom 4 with axiom $4'$. The definition of constant specification for $\LPCstandard$, and the definition of $\LPCstandard_\CS$ and  derivability in $\LPCstandard_\CS$, denoted by $\vdash_{\LPCstandard_\CS}$, is similar to that of $\LPC_\CS$. 

Relational models of $\LPCstandard$ are defined similar to $\LPC$-models satisfying a condition corresponding to axiom $4'$ (instead of condition 5 of Definition \ref{def: relational models conditions LPC}). More precisely, an $\LPCstandard_\CS$-model  $\model = (\worlds, \worldsNormal, \accrelFormula, \accrelTrem,  \valuation)$ is defined as an $\LPC_\CS$-model such that instead of condition 5 of Definition \ref{def: relational models conditions LPC} it satisfies the following condition:
\begin{description}
	\item[$5'.$] for all $w,v,u \in \worldsNormal$ and $u' \in \worlds$, for all $\phi, \psi \in \Formulae$, for all $s,t \in \Terms$: if $w \accrel_{\jbox{s}(\phi \limpliesCounterfactual \psi)} v$, $v \accrel_{\jbox{t}\phi} u$, and $u \accrel_{s \cdot t} u'$, then $u' \in \truthsetModel{\psi}$.
\end{description}
Validity is defined similar to Definition \ref{def: validity LPC} and is denoted by $\models_{\LPCstandard_\CS}$. 

\begin{theorem}[Completeness]
	Let $\CS$ be a constant specification for $\LPCstandard$. For each formula $\phi$ and set of formulas $T$,
	\[
	T \entails_{\LPCstandard_\CS}  \phi  \quad\text{iff}\quad  T \vdash_{\LPCstandard_\CS} \phi.
	\]
\end{theorem}
\begin{proof}
	The proof is similar to that given in Section \ref{sec: Completeness of LPC}. Again the major difference is that we have to take maximal $\LPCstandard_\CS$-consistent sets into consideration. The proof of the Truth Lemma is similar to that of Lemma \ref{lem: Truth Lemma LPC}. It is enough to show that the canonical relational model $\model = (\worlds, \worldsNormal, \accrelFormula, \accrelTrem,  \valuation)$ for $\LPCstandard_\CS$ is an $\LPCstandard_\CS$-model. Thus, we show that $\model$ satisfies the condition $5'$.
	
	Suppose that for $\Gamma, \Delta, \Sigma \in \worldsNormal$, $\Sigma' \in \worlds$,  $\phi, \psi \in \Formulae$, and $s, t \in \Terms$ we have that $\Gamma \accrel_{\jbox{s}(\phi \limpliesCounterfactual \psi)} \Delta$, $\Delta \accrel_{\jbox{t}\phi} \Sigma$, and $\Sigma \accrel_{s \cdot t} \Sigma'$. Since $\Gamma$ is a maximal $\LPCstandard_\CS$-consistent set, \linebreak $\jbox{s} (\phi \limpliesCounterfactual \psi) \limpliesCounterfactual (\jbox{t} \phi \limpliesCounterfactual \jbox{s \cdot t} \psi) \in \Gamma$. From $\Gamma \accrel_{\jbox{s}(\phi \limpliesCounterfactual \psi)} \Delta$
	it follows that $\jbox{t} \phi \limpliesCounterfactual \jbox{s \cdot t} \psi \in \Delta$. From this and $\Delta \accrel_{\jbox{t}\phi} \Sigma$, it follows that $\jbox{s \cdot t} \psi \in \Sigma$. From this and $\Sigma \accrel_{s \cdot t} \Sigma'$, it is deduced that $\psi \in \Sigma'$. By the Truth Lemma, it is immediate that $\Sigma' \in \truthsetModel{\psi}$, as desired. \qed  
\end{proof}

Similar to the extension $\LPCint$ of $\LPC$ in Section \ref{sec: The logic LPCint}, it is straightforward to extend the logic $\LPCstandard$ by axiom 9 to obtain a conditional justification logic that satisfies the internalization property. We leave the details to the reader.
\subsection{The logic $\LPCK$}
Normal conditional logics usually have the following rule in their axiomatization (see \cite{Chellas1975}):
\[
\lrule{ \phi \liff \psi}{ (\phi \limpliesCounterfactual \chi) \liff (\psi \limpliesCounterfactual \chi)}\, \RCEArule
\]
However, as the completeness proof of Section \ref{sec: Completeness of LPC} shows, we do not need the rule $\RCEArule$ in the axiomatization of $\LPC$. On the other hand, the following lemma shows that the rule $\RCEArule$ is not admissible in $\LPC$, and thus adding $\RCEArule$ to the rules of $\LPC$ gives a new logic.

\begin{lemma}
	The  rule
	\[
	\lrule{ \phi \liff \psi}{ (\phi \limpliesCounterfactual \chi) \liff (\psi \limpliesCounterfactual \chi)}\, \RCEArule
	\]
	is not validity preserving in $\LPC_\CS$, i.e. there exists $\phi, \psi \in \Formulae$ and $t \in \Terms$ such that $\entails_{\LPC_\CS} \phi \liff \psi$ but $\not \entails_{\LPC_\CS} (\phi \limpliesCounterfactual \chi) \liff (\psi \limpliesCounterfactual \chi)$.
\end{lemma}
\begin{proof}
	Let $\CS$ be an arbitrary constant specification for $\LPC$. Consider the relational model $\model = (\worlds, \worldsNormal, \accrelFormula, \accrelTrem,  \valuation)$ as follows:
	\begin{itemize}
		\item $\worlds = \worldsNormal = \{ w, v\}$,
		
		\item $\accrel_{p}(w) = \lset{w, v}$, $\accrel_{p \wedge p}(w) = \lset{w}$,\\
		$\accrel_{p}(v) = \accrel_{p \wedge p}(v) = \lset{v}$,\\
		$\accrel_{\phi}(w) = \accrel_{\phi}(v) =\truthsetModel{\phi}$, for all $\phi \in \Formulae$ such that $\phi \neq p$ and $\phi \neq p \wedge p$,
		
		\item $\accrel_{t}(w) = \accrel_{t}(v) = \{ w, v\}$, for all $t \in \Terms$,
		
		\item $\valuation(w) = \{ p, q \}$, $\valuation(v) = \{ p \}$.
	\end{itemize}
	It is easy to show that $\model$ is an $\LPC_\CS$-model. 
	
	We observe that $(\model,w) \not \entails p \limpliesCounterfactual q$ and $(\model,w) \entails p \wedge p \limpliesCounterfactual q$, and thus $\not \entails_{\LPC_\CS} (p \limpliesCounterfactual q) \liff (p \wedge p \limpliesCounterfactual q)$, but $\entails_{\LPC_\CS} p \liff p \wedge p$.  \qed
\end{proof}

Let $\LPCK$ denote the logic obtained from $\LPC$ by adding rule $\RCEArule$. The logic $\LPCK$  is an extension of the normal conditional logic $\Logic{CK}$+ID+MP of Chellas \cite{Chellas1975}. The definition of constant specification for $\LPCK$, and the definition of $\LPCK_\CS$ and derivability in $\LPCK_\CS$, denoted by $\vdash_{\LPCK_\CS}$, is similar to that of $\LPC_\CS$.

Relational models of $\LPCK$ are defined similar to $\LPC$-models satisfying an additional condition corresponding to the rule $\RCEArule$. More precisely, an $\LPCK_\CS$-model \linebreak $\model = (\worlds, \worldsNormal, \accrelFormula, \accrelTrem,  \valuation)$ is  defined as an $\LPC_\CS$-model that satisfies the following extra condition:
\begin{description}
	\item[$9.$] for all $\phi, \psi \in \Formulae$: if $\truthsetModel{\phi} \cap \worldsNormal = \truthsetModel{\psi} \cap \worldsNormal$, then $\accrel_{\phi} = \accrel_{\psi}$.
\end{description}
Validity is defined similar to Definition \ref{def: validity LPC} and is denoted by $\models_{\LPCK_\CS}$.

\begin{theorem}[Completeness]
	Let $\CS$ be a constant specification for $\LPCK$. For each formula $\phi$ and set of formulas $T$,
	\[
	T \entails_{\LPCK_\CS}  \phi  \quad\text{iff}\quad  T \vdash_{\LPCK_\CS} \phi.
	\]
\end{theorem}
\begin{proof}
	For soundness, suppose that $\models_{\LPCK_\CS} \phi \liff \psi$. From this it follows that in  every $\LPCK_\CS$-model $\model = (\worlds, \worldsNormal, \accrelFormula, \accrelTrem,  \valuation)$ we have $\truthsetModel{\phi} \cap \worldsNormal = \truthsetModel{\psi} \cap \worldsNormal$. From this and condition 9 we get $\accrel_{\phi} = \accrel_{\psi}$. Hence, for every $w \in \worldsNormal$ we get $\accrel_{\phi} (w) \subseteq \truthsetModel{\chi}$ iff $\accrel_{\psi} (w) \subseteq \truthsetModel{\chi}$. Therefore, $\models_{\LPCK_\CS} (\phi \limpliesCounterfactual \chi) \liff (\psi \limpliesCounterfactual \chi)$.
	
	For completeness, it suffices to show that the canonical relational model \linebreak $\model = (\worlds, \worldsNormal, \accrelFormula, \accrelTrem,  \valuation)$ satisfies condition 9. The proof of Truth Lemma is similar to that of Lemma \ref{lem: Truth Lemma LPC}. Furthermore, note that  if $\truthsetModel{\phi} \cap \worldsNormal = \truthsetModel{\psi} \cap \worldsNormal$, for some arbitrary $\phi, \psi \in \Formulae$, then $\vdash_{\LPCK_\CS} \phi \equiv \psi$. For suppose that $\not \vdash_{\LPCK_\CS} \phi \equiv \psi$. Then $\neg (\phi \equiv \psi) \in \Gamma$, for some $\Gamma \in \MCS_{\LPCK_\CS}$. Hence, either $\phi, \neg \psi \in \Gamma$ or $\neg \phi, \psi \in \Gamma$. In either case, by the Truth Lemma, we get $\truthsetModel{\phi} \cap \worldsNormal \neq \truthsetModel{\psi} \cap \worldsNormal$.
	
	Now in order to show condition 9, suppose that $\truthsetModel{\phi} \cap \worldsNormal = \truthsetModel{\psi} \cap \worldsNormal$. Thus, $\vdash_{\LPCK_\CS} \phi \equiv \psi$, and hence $(\phi \limpliesCounterfactual \chi) \liff (\psi \limpliesCounterfactual \chi)$ belongs to every maximal $\LPCK_\CS$-consistent set. Next we show that $\accrel_{\phi} = \accrel_{\psi}$. Suppose that $\Gamma \accrel_{\phi} \Delta$, for $\Gamma, \Delta \in \MCS_{\LPCK_\CS}$, and $\chi \in \Gamma/ \psi$. Thus, $\psi \limpliesCounterfactual \chi \in \Gamma$, and hence $\phi \limpliesCounterfactual \chi \in \Gamma$. Thus, $\chi \in \Delta$, and hence $\Gamma/ \psi \subseteq \Delta$. Therefore, $\Gamma \accrel_{\psi} \Delta$, and hence $\accrel_{\phi} \subseteq \accrel_{\psi}$. The proof of the converse is similar. \qed
\end{proof}

It is worth noting that the following rule from \cite{Chellas1975} is derivable in $\LPCK$
	\[
	\lrule{ \phi \liff \psi}{ (\chi \limpliesCounterfactual \phi) \liff (\chi\limpliesCounterfactual \psi )}\, (\sf RCEC)
	\]
	The proof is simple and thus is omitted here.

\subsection{The logic $\JFourC$}

In this part, we consider a fragment of $\LPC$ that will be used in Section \ref{sec: Nozick's analysis of knowledge} to formalize the justified true belief and counterfactual analyses of knowledge. 

The language of this fragment is the same as $\LPC$. Let $\JFourC$ denote the fragment of $\LPC$ which is obtained from $\LPC$ by dropping axiom 8, $\jbox{t} \phi \limpliesCounterfactual \phi$. Axiom 8 states that if it were that the agent knows $\phi$ for reason $t$, then it would be that $\phi$ is true. It is obvious that in the case of belief (instead of knowledge) this principle does not hold, because the agent may have a false belief.  Thus, $\JFourC$ can be considered as a justification logic of belief, where $\jbox{t} \phi$ is read ``the agent believes $\phi$ for reason $t$.'' 

The definition of constant specification for $\JFourC$ and the definition of $\JFourC_\CS$ are similar to that of $\LPC_\CS$. Relational models of $\JFourC$ are defined similar to $\LPC$-models in which requirement 6 of Definition \ref{def: relational models conditions LPC} is dropped. Validity is defined similar to Definition \ref{def: validity LPC}. 

Completeness again goes by way of a canonical model argument similar to that is given in Section \ref{sec: Completeness of LPC}, so we will omit the proof here.

\begin{theorem}[Completeness]
	The logic $\JFourC_\CS$ is sound and complete with respect to its relational models.
\end{theorem}

The logics ${\sf J4C}^{int}$, ${\sf J4C}^{st}$, and ${\sf J4CK}$ are defined in the same manner as the logics $\LPCint$, $\LPCstandard$, and $\LPCK$ are defined in the previous sections. In addition, one can define other fragments of $\LPC$ (such as $\JC$ and $\JTC$) and various extensions of $\LPC$ (such as $\Logic{JT45C}^+$, $\Logic{JBC}^+$, $\Logic{JDC}^+$, etc.) similar to the classical justification logics (for the fragments and extensions of $\LP$ see \cite{Art-Fit-Book-2019,Kuz-Stu-Book-2019}). For instance, the logic $\JC$ is obtained from $\LPC$ by removing axioms 8 and 9, and modifying the definition of constant specification as Definition 2.5 from \cite{Art-Fit-Book-2019}[page 16]. Relational models of $\JC$ are defined similar to $\LPC$-models in which the requirements 6 and 7 of Definition \ref{def: relational models conditions LPC} are dropped.\footnote{In Sections \ref{sec: Syntax JRC} and \ref{sec: Relational semantics for JRC} we introduce a conditional logic that its justification fragment is even weaker than $\JC$.}

\section{Nozick's analysis of knowledge}
\label{sec: Nozick's analysis of knowledge}

Before formalizing Nozick's definition of knowledge, we briefly recall the notion of Aumann knowledge that is usually used in Computer Science, Economics, and Game Theory applications. We then extend the language of $\JFourC$ with this knowledge modality.

The language of modal logic allows us to describe basic facts along with the knowledge and beliefs of agents. In the epistemic logic context, the knowledge of agents is sometimes formalized within the modal logic $\Logic{S5}$, based on Kripke models endowed with equivalence relations (see e.g. \cite{FHMV95,Meyer-Hoek-1995,DEL-2007}) or on Aumann’s partition models \cite{Aumann-1999}. In our relational models (that include both normal and non-normal states), this can be simply defined as the universal modality $\lknows$, quantifying universally over all epistemically possible normal states. Thus, $\lknows \phi$ is true (at any state) iff $\phi$ is true at all normal states. In other words, given a relational model $\model = (\worlds, \worldsNormal, \accrelFormula, \accrelTrem,  \valuation)$ the truth condition of $\lknows \phi$ is defined as follows:
\begin{equation}\label{eq: truth condition of Box}
(\model, w) \entails \lknows \phi   \text{ iff } (\model, v) \entails \phi, \text{ for all states }  v \in \worldsNormal.
\end{equation}
This notion of knowledge corresponds to the strong Leibnizian sense of ``truth in all possible worlds''.
The reason that for modeling knowledge the accessibility relation is
taken to be an equivalence relation, can be understood as follows: the
agent, being in a certain state, considers a set of alternatives which
are all alternatives of each other and one of which is the actual state
(so the agent considers his true state as an alternative).

In this setting, $\lknows$ is interpreted as “absolutely certain, infallible knowledge”, defined as truth in
all the worlds that are consistent with the agent’s information. This is a very strong form of knowledge encompassing all epistemic possibilities.
In this respect, $\lknows \phi$ could be best described as the possession of hard information. The $\Logic{S5}$ principles\footnote{We recall the axioms and rules of $\Logic{S5}$ involving the $\lknows$ modality here: the closure of knowledge:\linebreak $\lknows (\phi \limplies \psi) \limplies (\lknows \phi \limplies \lknows \psi)$; factivity: $\lknows \phi \limplies \phi$; positive introspection: $\lknows \phi \limplies \lknows \lknows \phi$; negative introspection: $\neg \lknows \phi \limplies \lknows \neg \lknows \phi$; and the necessitation rule: from $\phi$ infer $\lknows \phi$.\label{footnote: S5 axioms}}
 require that agents be fully introspective and immune to error, and thus the knowledge that the modal logic $\Logic{S5}$ is represented (i.e. the Aumann knowledge) is sometimes referred to as  `infallible knowledge' or `hard information' (see \cite{Aumann-1999,Baltag-Smets-LOFT-2008,Benthem-JANCL-2007}  for a more detailed exposition). In addition, the notion of knowledge that the universal modality provides suffers from the well-known logical omniscience problem: the agent knows all valid propositions (see \cite{FHMV95}). In other words, what the $\Logic{S5}$ knowledge modality is formalized is the propositions knowable to agent but not those that is actually known by him. This is indeed an idealized concept of knowledge. Putting these two facts together, the formula $\lknows \phi$ can be read ``the agent is entitled to infallibly know that $\phi$ is true.''

This notion of infallible knowledge is not very widely accepted by epistemologists. However, having the $\Logic{S5}$ knowledge modality is useful. It is indeed a simple way to add a knowledge modality to the language of $\JFourC$ without affecting the definition of relational models. Moreover, if an agent has infallible knowledge of a proposition then the agent's belief cannot be accidentally true. The condition of non-accidentality (i.e. the condition that if the agent knows a proposition then the agent's belief is not accidentally true) is proposed by McGinn \cite{McGinn-1999} as an adequacy condition on the definition of knowledge. We shall come back to this condition in the Gettier and McGinn examples below.

The language of justification logic, and hence that of $\JFourC$, can express the well-known justified true belief account of knowledge (JTB) as follows:\footnote{Note that a full formalization of JTB needs quantification on justification $\JTB(\phi):=  (\exists \justVarOne) \jbox{\justVarOne} \phi$ (see \cite{Fit08APAL} for a quantified logic of evidence). Since the language of all justification logics considered in this paper is propositional and has no quantifiers, we opt for the stronger formula $\JTB(\phi, t)$.}
\[
\JTB(\phi, t) := \phi \wedge \jbox{t} \phi.
\]
The formula $\JTB(\phi, t)$ says that ``$\phi$ is true  and it is also believed and justified by the evidence $t$.'' The reason why the truth component of JTB, i.e. the truth of $\phi$ in the formula $\JTB(\phi, t)$, has been separately stated is the fact that in the logic $\JFourC$ the justification assertion $\jbox{t} \phi$ only represents the justified belief of the agent.\footnote{In contrast to $\JFourC$, JTB can be defined in the logic $\LPC$ simply by $\jbox{t} \phi$. The fact that ``$\phi$ is true'' is implicitly stated in the formula $\jbox{t} \phi$, because the formula $\jbox{t} \phi \limpliesCounterfactual \phi$ is valid in $\LPC$. One may further show the validity of $\phi \wedge \jbox{t} \phi \Leftrightarrow \jbox{t} \phi$, and hence $\phi \wedge \jbox{t} \phi \equiv \jbox{t} \phi$, in $\LPC$.}

In the rest of this section, we reason in the combined logic $\JFourC \oplus \Logic{S5}$. In this hybrid logic, both the infallible knowledge and the justified true belief account of knowledge can be expressed. Furthermore, this hybrid logic enables us to formalize the counterfactual conditions on knowledge mentioned below. 

For simplicity, let $\Logic{L}$ denote the hybrid logic $\JFourC \oplus \Logic{S5}$. The language of $\Logic{L}$ extends that of $\JFourC$ with an operator $\lknows$ creating new formulas $\lknows \phi$. Axiomatically $\Logic{L}$ extends the axiomatic system $\JFourC$ with the axioms and rules of $\Logic{S5}$ from Footnote \ref{footnote: S5 axioms}. The definition of constant specification for $\Logic{L}$, and the definition of $\Logic{L}_\CS$ and derivability in $\Logic{L}_\CS$, denoted by $\vdash_{\Logic{L}_\CS}$, is similar to that of $\JFourC_\CS$. Models of $\Logic{L}$ are the same as $\JFourC$-models where the truth condition of formulas involving $\lknows$ is defined by \eqref{eq: truth condition of Box}. The various notions of validity and consequence relation from Definition \ref{def: validity LPC} carry over directly to the logic $\Logic{L}$.

We now turn our attention to Nozick's account of knowledge and its formalization in $\Logic{L}$. Nozick (\cite{Nozick-1981}[p. 179]) defines knowing via a method as follows:

An agent $S$ knows, via the method (or way of believing) $t$, that $\phi$ iff
\begin{enumerate}
	\item $\phi$ is true,
	
	\item $S$ believes, via method or way of coming to believe $t$, that $\phi$,
	
	\item If $\phi$ were false, then $S$ would not believe, via $t$, that $\phi$,
	
	\item  If $\phi$ were true, then $S$ would believe, via $t$, that $\phi$.
\end{enumerate}

We will refer to condition (3) as the sensitivity condition and to condition (4) as the
adherence condition. Since in this paper only mono-agent logics are considered, we neglect the agent $S$ in the above definition and formalize it as follows:

\[
K(\phi,t) := \phi \wedge \jbox{t} \phi \wedge (\neg \phi \limpliesCounterfactual \neg \jbox{t} \phi) \wedge (\phi \limpliesCounterfactual \jbox{t} \phi).
\]
Here, $\jbox{t} \phi$ is read ``the agent believes that $\phi$ is true by means of the method $t$". The first and second conjuncts formalize the justified true belief, and the  third and fourth conjuncts formalize the sensitivity and adherence conditions respectively.

It is known that the above account of knowledge resolves the Gettier examples (\cite{Gettier1963}). We only formalize Gettier's second example which is more simple. 

\begin{example}[Gettier's second example]\label{example: Gettier's second example}
	Suppose that Smith possesses good evidence in favor of the proposition that ``Jones owns a Ford'', for example, he has seen Jones driving a Ford. We denote this piece of evidence by $\justVarOne$. Smith also has a friend, Brown, whose current location is not known by Smith. Nonetheless, on the basis of his accepting that Jones owns a Ford, he infers and believes each of these two disjunctive propositions:
	\begin{itemize}
		\item ``Either Jones owns a Ford, or Brown is in Barcelona (or both)''.
		\item ``Either Jones owns a Ford, or Brown is in Boston (or both)''.
		
	\end{itemize}
	Smith realizes that he has good evidence for the first disjunct, and he sees this evidence as thereby supporting the disjunction as a whole. Suppose that, unknown to Smith, Jones does not own a Ford but only borrowed one from a friend, and in reality Brown is in Barcelona. Thus the first disjunction turns out to be (accidentally) true, and Smith has good evidence justifying it. But since there is significant luck in how Smith's belief manages to combine being true with being justified, it seems that Smith's justified true belief is not knowledge.
	
	Let us denote the sentence ``Jones owns a Ford'' by $p$, and denote by $q$ and $r$ the sentences ``Brown is in Barcelona'' and ``Brown is in Boston'' respectively. The set of assumptions in this example is
	\[
	G = \lset{  \neg p, q, \jbox{\justVarOne} p}.
	\]
	Here $\jbox{\justVarOne} p$ is read ``Smith has true belief that $p$ justified by $\justVarOne$.'' It is easy to show that 
	\begin{equation}\label{eq: Gettier example-has justified true belief}
		G \vdash_{\Logic{L}_\CS} (p \vee q) \wedge \jbox{\justConsThree \tapp \justVarOne} (p \vee q),
	\end{equation}
	where $\CS$ contains $\jbox{\justConsThree} (p \limplies p \vee q)$. This means that $\JTB(p \vee q, \justConsThree \tapp \justVarOne)$ is derivable from $G$. Thus, Smith has justified true belief that Jones owns a Ford or Brown is in Barcelona. However, the sensitivity condition does not hold, i.e. it is not the case that if $p \vee q$ were false, then Smith would not believe via $\justConsThree \tapp \justVarOne$ that $p \vee q$. Because at the state that is similar to the actual state but in which Brown is in Boston, we have that $p \vee q$ is not true but Smith still believes that  $p \vee q$. Thus, $\neg (p \vee q) \limpliesCounterfactual \neg \jbox{\justConsThree \tapp \justVarOne} (p \vee q)$ is false, and hence $K(p \vee q, \justConsThree \tapp \justVarOne)$ is false.

	The situation described in this example can modeled more precisely by the following  $\Logic{L}$-model. Define the relational model $\model = (\worlds, \worldsNormal, \accrelFormula, \accrelTrem,  \valuation)$ as follows:
	\begin{itemize}
		\item $\worlds = \worldsNormal =  \{ w, v\}$, 
		
		\item $\accrel_{\phi}(w) = \accrel_{\phi}(v) = \truthsetModel{\phi}$, for all $\phi \in \Formulae$ 
		
		\item $\accrel_{t}(w) = \accrel_{t}(v) =   \emptyset$ for all $t \in \Terms  $,
		
		\item $\valuation(w) = \{ q \}$ and $\valuation(v) = \{ r \}$.
	\end{itemize}
	\begin{figure}
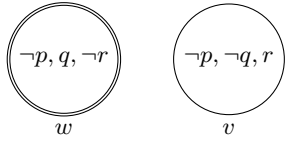

		\begin{center}
			\begin{tabular}{c}
				\begin{mytikz}
					
					\node[w,double,label={below:$w$}] (w) {$\neg p, q, \neg r$};
					
					\node[w,right of=w,label={below:$v$}] (v) {$\neg p, \neg q,r$};
					

				\end{mytikz}
			\end{tabular}
		\end{center}
		\caption{Model $\model$ for the Gettier's second example, in which $w$ is the actual state.}
	\end{figure}
	Among various possible states, we only consider two states in the model $\model$, at none of them Jones owns a Ford. At the state $w$, that is the actual state, Brown is in Barcelona, whereas at $v$, which is one of the most similar states to $w$, Brown is in Boston.
	
	It is easy to show that $\model$ is an $\Logic{L}_\CS$-model, where $\CS$ is any constant specification for $\Logic{L}$  (specially $\CS$ may contain $\jbox{\justConsThree} (p \limplies p \vee q)$).  We observe that $\accrel_{\justConsThree \tapp \justVarOne}(w) = \emptyset \subseteq \truthsetModel{p \vee q} = \lset{w}$ and thus 
	$$(\model,w) \entails (p \vee q) \wedge \jbox{\justConsThree \tapp \justVarOne} (p \vee q).$$
	This means that $(\model,w) \entails \JTB(p \vee q, \justConsThree \tapp \justVarOne)$. Note that all of the formulas of the set of assumptions $G$ are true at $w$.
	
	On the other hand, since $\accrel_{\neg (p \vee q)}(w) = \lset{v}$ and $\truthsetModel{\neg \jbox{\justConsThree \tapp \justVarOne} (p \vee q)} = \emptyset$, we have that $\accrel_{\neg (p \vee q)}(w)  \not \subseteq \truthsetModel{\neg \jbox{\justConsThree \tapp \justVarOne} (p \vee q)}$, and hence $(\model,w) \not \entails \neg (p \vee q) \limpliesCounterfactual \neg \jbox{\justConsThree \tapp \justVarOne} (p \vee q)$. Therefore, according to Nozick's account of knowledge, Smith does not know $\phi$. In other words, 
	\begin{equation}\label{eq: Gettier example-has no Nozick knowledge}
		\model, w \models \neg K(p \vee q,\justConsThree \tapp \justVarOne).
	\end{equation}
	Note that in this example we also have
	\begin{equation}\label{eq: Gettier example-has no infalliable knowledge}
		\model, w \models \neg \lknows (p \vee q).
	\end{equation}
	Thus Smith lacks both notions of knowledge. In fact, \eqref{eq: Gettier example-has justified true belief}, \eqref{eq: Gettier example-has no Nozick knowledge}, and \eqref{eq: Gettier example-has no infalliable knowledge} address the Gettier problem, by showing that Smith has justified true belief but has no knowledge, and this is compatible with our intuition. \qed
\end{example}

However, it is known that Nozick's account of knowledge is not necessary and sufficient for knowledge. Consider the following example from McGinn \cite{McGinn-1999}. He provides an example to show that Nozick's counterfactuals are not sufficient for knowledge, i.e. the sensitivity condition holds but it is not intuitive to attribute knowledge to the agent.

\begin{example}[McGinn's example]\label{example: McGinn's example}
	``You visit a hitherto unexplored country in which the inhabitants have the custom of simulating being in pain. You do not know that their pain behaviour is mere pretence, and so you form the belief of each person you meet that he or she is in pain; imagine you have	acquired a great many false beliefs in this way. There is, however, one person in this country who is an exception to the
	custom of pain pretence: this hapless individual is in constant pain and shows it (we can suppose that he falsely	believes others to be in his unfortunate condition—he has not been told of the pretence by the others). You also believe of this person, call him $N$, that he is in pain. Now I take it that we would not say that your true belief that $N$ is in pain counts as knowledge, for it is, intuitively, a mere accident that your belief is true in this instance. But now	consider the relevant counterfactuals, in particular ‘if $N$ were not in pain, you would not believe that $N$ was in pain’: this counterfactual is true in the envisaged circumstances, since if $N$ were not in pain then (unlike the pretenders around him) he would not behave as if he was, and so you would not believe that he was. So your belief that $N$ is in	pain does track the truth of that proposition even though it does not rank as knowledge.'' (Quoted from \cite{McGinn-1999}.)
	
	Let $p$ denote the proposition ``$N$ is in pain'' and $\justVarOne$ denote ``your observation of the fact that $N$ is in	pain.'' Here $\jbox{\justVarOne} p$ is read ``You have true belief that $p$, justified by reason $\justVarOne$.'' Define the relational model $\model' = (\worlds', \worldsNormal', \accrelFormula', \accrelTrem',  \valuation')$ as follows:
	\begin{itemize}
		\setlength\itemsep{0.1cm}
		\item $\worlds' = \worldsNormal' = \lset{w,v}$, 
		
		
		\item $\accrel'_{\phi}(w) = \accrel'_{\phi}(v) = \truthset{\phi}{{\model'}}$,  for all $\phi \in \Formulae$
		
		\item $\accrel'_{t}(w)  = \lset{w}$ and $ \accrel'_{t}(v) = \lset{v}$ for all $t \in \Terms  $,
		
		\item $\valuation'(w) = \{ p \}$, $\valuation'(v) =\emptyset$.
	\end{itemize}
	It is easy to show that $\model'$ is an $\Logic{L}_\CS$-model, where $\CS$ is any constant specification for $\Logic{L}$. Indeed, at the model $\model'$ the state $v$ is essentially the same as the actual state $w$ except that $\neg p$ is true at $v$, i.e. $N$ is not in pain at $v$.
	
	We observe that $(\model',v) \not \entails p$ and thus $(\model',w) \not \entails \lknows p$. On the other hand, since \linebreak $\accrel'_{\neg p}(w) =  \lset{v}$ and $\truthset{\neg \jbox{\justVarOne} p}{{\model'}} = \lset{v}$, we have that $\accrel'_{\neg p}(w)   \subseteq \truthset{\neg \jbox{\justVarOne} p}{{\model'}}$, and hence \linebreak $(\model',w)  \entails \neg p \limpliesCounterfactual \neg \jbox{\justVarOne} p$. 
	
	In addition, since  $\accrel'_{p}(w) =  \lset{w}$ and $\truthset{\jbox{\justVarOne} p}{{\model'}} = \lset{w}$, we have that $\accrel'_{p}(w)   \subseteq \truthset{\jbox{\justVarOne} p}{{\model'}}$, and hence $(\model',w)  \entails  p \limpliesCounterfactual  \jbox{\justVarOne} p$. 

	Thus, although both the sensitivity condition and the adherence condition hold (at $w$) the Aumann knowledge does not hold. Therefore, neither the sensitivity condition nor the adherence condition are sufficient for Aumann knowledge.
	
	As we have already observed $(\model',w) \not \entails \lknows p$. It is also not difficult to show that \linebreak $(\model',w)  \entails K(p,\justVarOne)$. Putting these two facts together, it follows that Nozick knowledge does not imply Aumann knowledge. \qed
\end{example}

Next, we show that the Aumann knowledge does not imply Nozick's account of knowledge as well.

\begin{example}\label{example: infallible knowledge does not imply Nozick knowledge}
Define the relational model $\model'' = (\worlds'', \worldsNormal'', \accrelFormula'', \accrelTrem'',  \valuation'')$ as follows:
\begin{itemize}
	\setlength\itemsep{0.1cm}
	\item $\worlds'' = \{ w, v\}$; $\worldsNormal'' = \{ w \}$,
	
	\item $\accrel''_{\phi}(w) = \truthsetModel{\phi} \cap \worldsNormal$, for all $\phi \in \Formulae$,
	
	\item $\accrel''_{\justVarOne}(w) = \{ w, v\}$, for all $\justVarOne \in \VTerms$; 
	$\accrel''_{\justConsThree}(w) = \{ w \}$  for all $\justConsThree \in \CTerms$;
	$\accrel''_{t}(w) = \emptyset$, for all $t \in \Terms \setminus (\VTerms \cup \CTerms)$; 
	$\accrel''_{t}(v) = \emptyset$, for all $t \in \Terms$,
	
	\item $\valuation''(w) = \{ p \}$, $\valuation''(v) = \emptyset$.

\end{itemize}
It is easy to show that $\model''$ is an $\Logic{L}_\CS$-model, where $\CS$ is any constant specification for $\Logic{L}$. We observe that $(\model'',w) \entails \lknows p$. On the other hand, since $\accrel''_{ p}(w) =  \truthset{p}{{\model''}}\cap \worldsNormal = \lset{w}$ and $\truthset{\jbox{\justVarOne} p}{{\model''}} = \emptyset$, we have that $\accrel''_{ p}(w)  \not \subseteq \truthset{ \jbox{\justVarOne} p}{{\model''}}$, and hence $(\model'',w) \not \entails  p \limpliesCounterfactual  \jbox{\justVarOne} p$. Therefore, $(\model'',w) \not \entails K(p,\justVarOne)$.
\qed
\end{example}

Let us summarize the results. Let $\CS$ be an arbitrary constant specification for $\Logic{L}$. The relational models $\model$, $\model'$, and $\model''$ of Examples \ref{example: Gettier's second example}--\ref{example: infallible knowledge does not imply Nozick knowledge} show the following facts:
\begin{itemize}
	\item The JTB account of knowledge is not knowledge:
	\begin{eqnarray*}
		\JTB(\phi, t) &\not \entails_{\Logic{L}_\CS}& \lknows \phi, \qquad~~~~~ \text{ for some $\phi$ and $t$,}\\
		\JTB(\phi, t) &\not \entails_{\Logic{L}_\CS}& K(\phi, t), \qquad \text{ for some $\phi$ and $t$.}
	\end{eqnarray*}
	
	\item Neither the sensitivity condition nor the adherence condition is sufficient for Aumann knowledge:
	
	\begin{eqnarray*}
		(\neg \phi \limpliesCounterfactual \neg \jbox{t} \phi) &\not \entails_{\Logic{L}_\CS}& \lknows \phi, \qquad \text{ for some $\phi$ and $t$,}\\
		(\phi \limpliesCounterfactual \jbox{t} \phi) &\not \entails_{\Logic{L}_\CS}& \lknows \phi, \qquad \text{ for some $\phi$ and $t$.}
	\end{eqnarray*}
	
	\item Nozick knowledge does not imply Aumann knowledge and vice versa:
	\begin{eqnarray*}
		K(\phi, t) &\not \entails_{\Logic{L}_\CS}& \lknows \phi, \qquad~~~~~ \text{ for some $\phi$ and $t$,}\\
		\lknows \phi &\not \entails_{\Logic{L}_\CS}& K(\phi, t), \qquad \text{ for some $\phi$ and $t$.}
	\end{eqnarray*}
	
\end{itemize}


\section{Justification logic with a relevant counterfactual conditional}
\label{sec: Syntax JRC}

In this section, we introduce a justification logic with a relevant counterfactual conditional. This conditional is inspired by the work of Priest in \cite{Priest2008}[pages 208--211]. By adding the Routley star to relational models and by applying some modifications to $\LPC$-models, we introduce a justification logic, called $\JRC$, with a relevant counterfactual, denoted here by $\limpliesCounterRel$. The language of $\JRC$ also includes a relevant conditional, denoted by $\limpliesRelevant$. The logic $\JRC$ is an extension of the relevant logic $\Logic{B}$ and the conditional logic $\Logic{C}^+$.\footnote{It is worth remarking that Pereira and Apar\'{\i}cio \cite{Pereira-Aparicio-1989} define a relevant counterfactual $\limpliesCounterfactual_R$ as\linebreak $\phi \limpliesCounterfactual_R \psi := (\phi \limpliesCounterfactual \psi) \wedge \neg (\phi \limpliesCounterfactual (\neg \psi \limpliesCounterfactual(\phi \limpliesCounterfactual \neg \psi)))$. Later, Mares and Fuhrmann \cite{Mares-Topoi-1994,Mares-Fuhrmann-1995} define a relevant counterfactual, here denoted by  $\limpliesCounterfactual_{MF}$, as $\phi \limpliesCounterfactual_{MF} \psi := \phi \limpliesCounterfactual (\phi \limpliesRelevant \psi)$. One may extend justification logics with either $\limpliesCounterfactual_R$ or $\limpliesCounterfactual_{MF}$. We leave the study of these extensions to possible future work. }

Justification logics are usually formulated as axiom systems. However, $\JRC$ is primarily formulated as a tableau proof system. For this reason, there is no need to relativize $\JRC$ by constant specifications. Moreover, in contrast to the conditional justification logics of the previous sections, the language of $\JRC$ does not contain the application operation `$\tapp$', since there is no need to formalize Modus Ponens in the language by the application axiom.

Let us now go over the formal definition of the syntax of $\JRC$. Justification terms are constructed from a countable set of justification variables $\VTerms$ and formulas are constructed from a countable set of atomic propositions $\Prop$ by the following grammars:  
\begin{gather*}
	t \coloncolonequals   \justVarOne \in \VTerms \mid  t + t \, , \\
	\phi \coloncolonequals p \in \Prop \mid\ \lnegRelevant \phi \mid \phi \wedge \phi \mid \phi \limpliesRelevant \phi \mid \phi \limpliesCounterRel \phi  \mid \jbox{t} \phi \, .
\end{gather*}
The set of justification terms and formulas are denoted by  $\TermsJCR$ and $\FormulaeJCR$ respectively. Here, $\lnegRelevant$ denotes the intensional negation of relevant logics, and $\limpliesCounterRel$ denotes the relevant counterfactual, and $\limpliesRelevant$ denotes the relevant conditional.

The extensional disjunction is defined by  $\phi \vee \psi \colonequals \lnegRelevant (\lnegRelevant \phi \wedge \lnegRelevant \psi)$, and fusion (or intensional conjunction) is defined by $\phi \lfussion \psi \colonequals \lnegRelevant (\phi \limpliesRelevant \lnegRelevant \psi)$.

For a  formula $\phi$, the
set of all subformulas of $\phi$, denoted by $Sub(\phi)$, is defined inductively as follows:
$ Sub(p)=\{p\}$, for any atomic proposition $p$;
$Sub(\lnegRelevant \phi)=\{\lnegRelevant \phi\}\cup Sub(\phi)$; $Sub(\phi \star \psi)=\{\phi \star \psi\}\cup Sub(\phi)\cup Sub(\psi)$ (where $\star$ is either $\wedge$, $\limpliesRelevant$, or $\limpliesCounterRel$); $Sub(\jbox{t} \phi)=\{\jbox{t} \phi\}\cup Sub(\phi)$. 

\section{Relational semantics for $\JRC$}
\label{sec: Relational semantics for JRC}

The semantics of $\JRC$ is defined similar to that of $\LPC$. The major difference is that the Ruotley star is added to define the truth condition of the intentional negation $\lnegRelevant$, a ternary relation is added to define the truth condition of the relevant conditional $\limpliesRelevant$, and some changes have been made in the definitions of $\accrelFormula$ and valuation function. For the relevant conditional, we employ a simplified version of the Routly-Meyer ternary relation semantics given in \cite{Priest-Sylvan-1992,Priest2008}.

\begin{definition}\label{def: Routley Relational models}
	A Routley relational model is a tuple $\model= (\worlds, \worldsNormal, \relRelevant, \accrelFormula, \accrelTrem, \ast, \valuationR)$ where
	
	\begin{enumerate}
		\item $W$ and $W_N$ are non-empty sets of states and normal states, respectively, such that $W_N \subseteq W$.

		\item $\relRelevant \subseteq \worlds \times \worlds \times \worlds$ is a ternary relation satisfying the \textit{normality condition}: 
		\[
		\mbox{if $w \in \worldsNormal$, then } (\relRelevant wvu \mbox{ iff } v=u).
		\]
		
		\item $\accrelFormula$ is  a function that assigns to each formula $\phi \in \FormulaeJCR$ a binary relation $\relFormula$ on $\worlds$. Given $w \in \worlds$, let $\relFormula(w) := \{ v \in \worlds \mid w \relFormula v \}$.

		\item $\accrelTrem$ is  a function that assigns to each term $t \in \TermsJCR$ a binary relation $\relTrem$ on $\worlds$. Given $w \in \worlds$, let $\relTrem(w) := \{ v \in \worlds \mid w \relTrem v \}$. 
		
		\item $*: \worlds \to \worlds$ is a function such that for all $w \in \worlds$
		\[
		w^{**} = w.
		\]

		\item $\valuationR$ is a valuation function such that $\valuationR: \Prop \to \powerset(\worlds)$. 
	\end{enumerate}
\end{definition}

\begin{definition}
	Let $\model= (\worlds, \worldsNormal, \relRelevant, \accrelFormula, \accrelTrem, \ast, \valuationR)$ be a Routley  relational model. Truth at states $w \in \worlds$ are inductively defined as follows:
	\begin{align*}
		(\model, w) &\entails p \text{ iff } w \in \valuationR(p) \, ,\\
		(\model, w) &\entails \lnegRelevant \phi \text{ iff } (\model, w^\ast) \not\entails \phi \, ,\\
		(\model, w) &\entails \phi \wedge \psi \text{ iff } (\model, w) \entails \phi \text{ and } (\model, w) \entails \psi \, ,\\
		(\model, w) &\entails \phi \limpliesRelevant \psi \text{ iff for all } v,u \in \worlds \text{ such that } \relRelevant wvu, \text{ if } (\model, v) \entails \phi, \text{ then } (\model, u) \entails \psi  \, ,\\
		(\model, w) &\entails \phi \limpliesCounterRel \psi \text{ iff }  \text{ for all  }  v \in \worlds,  \text{if }   w \relFormula  v, \text{ then }  (\model, v) \entails \phi\, ,\\
		(\model, w) &\entails \jbox{t} \phi \text{ iff }  \text{ for all  }  v \in \worlds, \text{ if }   w \relTrem  v, \text{ then }  (\model, v) \entails \phi.		
	\end{align*}
	
\end{definition}

Comparing this definition with Definition \ref{def: Truth LPC}, we see that the laws of logic hold in non-normal states of Routley relational models. As we shall see below, normal states will be only used in the conditions imposed on the models (see Definition \ref{def: relational models conditions JRC}) and in the definition of validity (see Definition \ref{def: validity JRC}).

\begin{definition}\label{def: relational models conditions JRC}
	A Routley relational model $\model= (\worlds, \worldsNormal,\relRelevant, \accrelFormula, \accrelTrem, \ast, \valuationR)$ is called a $\JRC$-model if it satisfies the following conditions:
	\begin{enumerate}
		\item For all $w \in \worldsNormal$, $\relFormula(w) \subseteq \truthsetModel{\phi}$.
		
		\item For all $w \in \worlds$, if $w \in \truthsetModel{\phi}$, then $w \in \relFormula(w)$.
		
		\item For all $s,t \in \TermsJCR$, $\accrel_{s \ttsum t} \subseteq \accrel_{s} \cap \accrel_{t}$.

	\end{enumerate}
\end{definition}

\begin{definition}\label{def: validity JRC}
	
	\begin{enumerate}
		\item Given a $\JRC$-model $\model= (\worlds, \worldsNormal, \relRelevant, \accrelFormula, \accrelTrem, \ast, \valuationR)$, we write $\model \entails \phi$ if
		for all $w \in \worldsNormal$, we have 
		$(\model, w) \entails \phi$.
		
		\item A formula $\phi$ is $\JRC$-valid, denoted by $\entails_{\JRC} \phi$, if $\model \entails \phi$ for all $\JRC$-models $\model$.
		
		\item Given a set of formulas $T$ and a formula $\phi$ of $\JRC$, the (local) consequence relation is defined as follows: \\ $T \entails_{\JRC} \phi$ iff for all 
		$\JRC$-models $\model= (\worlds, \worldsNormal,\relRelevant, \accrelFormula, \accrelTrem, \ast, \valuationR)$, for all $w \in \worldsNormal$, if $(\model, w) \entails \psi$ for all $\psi \in T$, then $(\model, w) \entails \phi$.
	\end{enumerate}
	
\end{definition}

It can be shown that all axioms of the basic relevant logic $\Logic{B}$ are $\JRC$-valid and all rules of $\Logic{B}$ preserve the $\JRC$-validity. Thus, $\JRC$ is an extension of $\Logic{B}$. 

In the next proposition, we list some valid sequents involving the relevant counterfactual $\limpliesCounterRel$. We omit its rather simple proof.

\begin{proposition}

	\begin{enumerate}

		\item $\entails_{\JRC} \phi \limpliesCounterRel \phi$.
		
		\item $\entails_{\JRC} \phi \wedge \psi \limpliesCounterRel \phi$, $\entails_{\JRC} \phi \wedge \psi \limpliesCounterRel \psi$.
		
		\item $\entails_{\JRC} ((\phi \limpliesCounterRel \psi) \wedge (\phi \limpliesCounterRel \chi)) \limpliesCounterRel (\phi \limpliesCounterRel (\psi \wedge \chi))$.
		
		\item $\entails_{\JRC}\ \lnegRelevant \lnegRelevant \phi \limpliesCounterRel \phi$, $\entails_{\JRC}  \phi \limpliesCounterRel\ \lnegRelevant \lnegRelevant \phi$.
		
		\item $\phi, \phi \limpliesCounterRel \psi \models_{\JRC} \psi$.

		\item  $\entails_{\JRC} \jbox{s} \phi \limpliesCounterRel  \jbox{s \ttsum t} \phi$, $\entails_{\JRC} \jbox{t} \phi \limpliesCounterRel  \jbox{s \ttsum t} \phi$.

	\end{enumerate}
\end{proposition}

It is worth noting that none of the paradoxes of material implication, including $\phi \limpliesCounterRel (\psi \limpliesCounterRel \phi)$ and $\lnegRelevant \phi \limpliesCounterRel (\phi \limpliesCounterRel \psi)$, are valid in $\JRC$.

\section{Relevance of $\limpliesCounterRel$}

It is known that if $\phi \limpliesRelevant \psi$ is valid in a relevant logic, then $\phi$ and $\psi$ have at least one atomic proposition in common (see e.g. \cite{Avron-APAL-2014}). This is called Belnap's variable-sharing property (see \cite{Anderson-Belnap-1975}). We show that the same property holds for the conditional $\limpliesCounterRel$.

\begin{theorem}
	If $\entails_{\JRC}  \phi \limpliesCounterRel \psi$, then $\phi$ and $\psi$ have at least one atomic proposition in common.
\end{theorem}

\begin{proof}
	Suppose that $\entails_{\JRC}  \phi \limpliesCounterRel \psi$, and assume to obtain a contradiction that  $\prop{\phi} \cap \prop{\psi} = \emptyset$. Define the model $\model= (\worlds, \worldsNormal, \relRelevant, \accrelFormula, \accrelTrem, \ast, \valuationR)$ as follows:
	\begin{itemize}
		\item Let $\worlds = \lset{ w, v, v^\ast }$ and $\worldsNormal = \lset{ w }$, and let $w^\ast = w$.
		
		\item Let $\relRelevant = \lset{ (w,w,w), (w,v,v), (w, v^*, v^*), (v,v,v), (v,v^*,v), (v^*,v,v^*)}$.
		
		\item For all $\chi \in \FormulaeJCR$, let $\accrel_{\chi} (w) = \truthsetModel{ \chi }$, $\accrel_{\chi} (v) = \lset{ v }$, and $\accrel_{\chi} (v^\ast) = \lset{ v^\ast }$.
		
		\item For all $t \in \TermsJCR$, let $\accrel_{t} (w) = \emptyset$, $\accrel_{t} (v) = \lset{ v }$, and $\accrel_{t} (v^\ast) = \lset{ v^\ast }$.
		
		\item For all $p \in \prop{\phi}$, let $\valuationR(p) = \lset{ v }$, and for all $p \in \prop{\psi}$, let $\valuationR(p) = \lset{ v^\ast }$. For other atomic propositions $p$, the definition of $\valuationR(p)$ is arbitrary.
	\end{itemize}
	We can prove the following facts:
	\begin{description}
		\item[Fact 1.] $\Sub{\phi} \cap \Sub{\psi} = \emptyset$.
		
		\item[Fact 2.] $\model$ is a $\JRC$-model.
		
		\item[Fact 3.] For all $\phi' \in \Sub{\phi}$, we have $(\model, v) \models \phi'$ and $(\model, v^\ast) \not \models \phi'$.
		
		\item[Fact 4.] For all $\psi' \in \Sub{\psi}$, we have $(\model, v) \not \models \psi'$ and $(\model, v^\ast)  \models \psi'$.
	\end{description}
	The proof of Fact 1 is trivial. Let us show Fact 2.  
	We now show that conditions 1--3 of Definition \ref{def: relational models conditions JRC} hold in $\model$. 
	
	From  $\accrel_{\chi} (w) = \truthsetModel{ \chi }$, it follows that condition 1 holds. Again since  $\accrel_{\chi} (w) = \truthsetModel{ \chi }$, $v \in \accrel_{\chi} (v)$, and $v^\ast \in \accrel_{\chi} (v^\ast)$ it follows that condition 2  holds. Since for all $t \in \TermsJCR$ we have $\accrel_{t} = \lset{ (v,v), (v^\ast, v^\ast) }$, condition 3  clearly holds. 
	
	Fact 3 is proved by induction on the complexity of $\phi'$.
	
	The case where $\phi'$ is an atomic proposition follows from the definition of $\valuationR$. For the induction step, we consider in detail only three cases (other cases are similar). 
	
	Suppose that $\phi' = \phi_1 \limpliesRelevant \phi_2$. Since $\phi_1, \phi_2 \in \Sub{\phi}$, by the induction hypothesis we have $(\model, v) \models \phi_1$, $(\model, v) \models \phi_2$, $(\model, v^\ast) \not \models \phi_1$, and $(\model, v^\ast) \not \models \phi_2$. From these facts and  $\relRelevant vvv$, $\relRelevant vv^*v$, and $\relRelevant v^*vv^*$, it follows that  $(\model, v) \models \phi_1 \limpliesCounterRel \phi_2$ and  $(\model, v^\ast) \not\models \phi_1 \limpliesCounterRel \phi_2$.
	
	Suppose that $\phi' = \phi_1 \limpliesCounterRel \phi_2$. Since $\phi_1, \phi_2 \in \Sub{\phi}$, by the induction hypothesis we have $(\model, v) \models \phi_1$, $(\model, v) \models \phi_2$, $(\model, v^\ast) \not \models \phi_1$, and $(\model, v^\ast) \not \models \phi_2$. Since $\accrel_{\phi_1} (v) = \lset{ v }$ and $\truthsetModel{\phi_2}$ is either $\lset{ v }$ or $\lset{ w, v }$, it follows that $\accrel_{\phi_1} (v) \subseteq \truthsetModel{\phi_2}$, and hence $(\model, v) \models \phi_1 \limpliesCounterRel \phi_2$. On the other hand, Since $\accrel_{\phi_1} (v^\ast) = \lset{ v^\ast }$ and $\truthsetModel{\phi_2}$ is either $\lset{ v }$ or $\lset{ w, v }$, it follows that $\accrel_{\phi_1} (v^\ast) \not\subseteq \truthsetModel{\phi_2}$, and hence $(\model, v^\ast) \not\models \phi_1 \limpliesCounterRel \phi_2$.
	
	Suppose that $\phi' = \jbox{t} \chi$. Since $\chi \in \Sub{\phi}$, by the induction hypothesis we have $(\model, v) \models \chi$  and $(\model, v^\ast) \not \models \chi$. Since $\accrel_{t} (v) = \lset{ v }$ and $\truthsetModel{\chi}$ is either $\lset{ v }$ or $\lset{ w, v }$, it follows that $\accrel_{t} (v) \subseteq \truthsetModel{\chi}$, and hence $(\model, v) \models \jbox{t} \chi$. On the other hand, Since $\accrel_{t} (v^\ast) = \lset{ v^\ast }$ and $\truthsetModel{\phi_2}$ is either $\lset{ v }$ or $\lset{ w, v }$, it follows that $\accrel_{t} (v^\ast) \not\subseteq \truthsetModel{\chi}$, and hence $(\model, v^\ast) \not\models \jbox{t} \chi$.
	
	Fact 4 is proved by induction on the complexity of $\psi'$. The proof is quite similar to the proof of Fact 3 and can be safely omitted here.
	
	Finally, a contradiction follows from Facts 1--4. Since $\phi \in \Sub{\phi}$ and $\psi \in \Sub{\psi}$, from Facts 3 and 4 we get $(\model, v) \models \phi$ and $(\model, v) \not\models \psi$.	Thus, $v \in \truthsetModel{\phi}$ but $v \not \in \truthsetModel{\psi}$. Since $\accrel_{\phi} (w) = \truthsetModel{ \phi }$, we conclude that  $\accrel_{\phi} (w) \not \subseteq \truthsetModel{ \psi }$. Therefore, $(\model, w)\not \models \phi \limpliesCounterRel \psi$, which contradicts the assumption $\entails_{\JRC}  \phi \limpliesCounterRel \psi$.
	\qed
\end{proof}
\section{Tableaux for $\JRC$}
\label{sec: Tableaux for JRC}

As usual, a tableau for a formula is a rooted binary tree constructed by applying the tableau rules. Tableaux of $\JRC$ contain expressions of the following forms
\[
\text{ $\phi, +x$ or $\phi, -x$, or $x \relTableau{\phi} y$ or $x \relTableau{t} y$ or $rxyz$}
\]
where $x, y, z \in \lset{ 0, 1, 2, \ldots } \cup \lset{ 0^\sharp, 1^\sharp, 2^\sharp, \ldots}$. Moreover, let
\[
\bar{x} = \left\{
\begin{array}{ll}
	i & \mbox{ if $x = i^\sharp$ } \\
	i^\sharp & \mbox{ if $x = i$ }
\end{array}
\right.
\]

\begin{figure}
	\centering
	\begin{tabular}{|c |c |c |c |c |c |c |}
		\hline
		$(T \sim)$ & $(T \wedge)$ & $(T \limpliesRelevant)$& $(T \limpliesCounterRel)$ & $(Cut_r)$ & $(T\jbox{}\,)$ & $(\relTableau{+})$  \\
		\hline
		\alwaysNoLine\AXC{}
		\UIC{$\lnegRelevant \phi, +x$}
		\UIC{$|$}
		\UIC{$\phi, - \bar{x}$}\DP
		
		&
		\alwaysNoLine
		\AXC{}
		\UIC{$\phi\wedge \psi, +x$}
		\UIC{$|$}
		\UIC{$\phi, +x$}
		\UIC{$ \psi, +x$}
		\UIC{}\DP
		&
		\alwaysNoLine\AXC{}
		\UIC{$\phi \limpliesRelevant \psi, +x$}
		\UIC{$rxyz$}
		\UIC{$\diagup\diagdown$}
		\UIC{$\phi, -y~~~ \psi, +z$}\DP
		&
		\alwaysNoLine\AXC{}
		\UIC{$\phi \limpliesCounterRel \psi, +x$}
		\UIC{$x \relTableau{\phi} y$}
		\UIC{$|$}
		\UIC{$ \psi, +y$}\DP
		&
		\alwaysNoLine\AXC{}
		\UIC{$$}
		\UIC{$\diagup\diagdown$}
		\UIC{$\phi, -x~~~~ \phi, +x$}
		\UIC{$~~~~~~~~~~~~~~ x \relTableau{\phi} x$}
		\DP
		&
		\alwaysNoLine\AXC{}
		\UIC{$\jbox{t} \phi, +x$}
		\UIC{$x \relTableau{t} y$}
		\UIC{$|$}
		\UIC{$ \phi, +y$}
		\DP
		&
		\alwaysNoLine\AXC{}
		\UIC{$x \relTableau{s \ttsum t} y$}
		\UIC{$|$}
		\UIC{$x \relTableau{s} y$}
		\UIC{$x \relTableau{t} y$}\DP 
		\\
		\hline \hline
		$(F \sim)$ & $(F \wedge)$ & $(F \limpliesRelevant)$& $(F \limpliesCounterRel)$ & $(F \limpliesCounterRel_0)$ & $(F\jbox{}\,)$ & \text{Normality rule}\\
		\hline
		\alwaysNoLine\AXC{}
		\UIC{$\lnegRelevant \phi, -x$}
		\UIC{$|$}
		\UIC{$\phi, + \bar{x}$}\DP 
		
		&
		\alwaysNoLine\AXC{}
		\UIC{$\phi\wedge \psi, -x$}
		\UIC{$\diagup\diagdown$}
		\UIC{$\phi, -x~~~ \psi, -x$}\DP
		&
		\alwaysNoLine
		\AXC{}
		\UIC{$\phi \limpliesRelevant \psi, -x$}
		\UIC{$|$}
		\UIC{$rxjk$}
		\UIC{$\phi, +j$}
		\UIC{$\psi, -k$}
		\DP
		&
		\alwaysNoLine
		\AXC{}
		\UIC{$\phi \limpliesCounterRel \psi, -x$}
		\UIC{$|$}
		\UIC{$x \relTableau{\phi} j$}
		\UIC{$\psi, -j$}
		\DP
		&
		\alwaysNoLine
		\AXC{}
		\UIC{$\phi \limpliesCounterRel \psi, -0$}
		\UIC{$|$}
		\UIC{$0 \relTableau{\phi} j$}
		\UIC{$ \phi, +j$}
		\UIC{$\psi, -j$}
		\UIC{}
		\DP
		&
		\alwaysNoLine
		\AXC{}
		\UIC{$\jbox{t} \phi, -x$}
		\UIC{$|$}
		\UIC{$x \relTableau{t} j$}
		\UIC{$\phi, -j$}
		\DP
		&
		\alwaysNoLine
		\AXC{}
		\UIC{$\cdot$}
		\UIC{$|$}
		\UIC{$r0xx$}
		\UIC{}
		\DP
		\\\hline
	\end{tabular}
	\caption{Tableau rules of $\JRC$. In the rules, $x, y, z$ are anything of the form  $i$ or $i^\sharp$, for some non-negative integer $i$. And $j, k$ are non-negative integers.} \label{Figure: Tableau rules of JRC}
\end{figure}

Tableau rules of $\JRC$ are presented in Figure \ref{Figure: Tableau rules of JRC}. All the rules, except the rules $(T\jbox{})$, $(F\jbox{ })$, and $(\relTableau{+})$, comes from \cite{Priest2008}. It is worth noting that the rule $(Cut_r)$ is formulated only for $x=0$ in \cite{Priest2008}[page 210]. 

\noindent
\textbf{Rule restrictions.} The tableau rules of $\JRC$ should meet the following requirements:

\begin{itemize}
	\item In the rule $(F \limpliesRelevant)$, we have that  $j, k \in \mathbb{N}$ are new and if $x = 0$ then $j=k$.
	
	\item In the rule $(F \limpliesCounterRel)$, we have that $x \neq 0$ and $j \in \mathbb{N}$ is new.
	
	\item In the rules $(F \limpliesCounterRel_0)$ and $(F\jbox{}\,)$, we have that $j \in \mathbb{N}$ is new.
	
	\item In the rule $(Cut_r)$, $\phi$ occurs as the antecedent of a conditional $\limpliesCounterRel$ and $x$ occurs on the branch.
	
	\item In the normality rule, $x$ occurs on the branch.
\end{itemize}

A tableau branch closes if a pair of the form $\phi, +x$, and $\phi, -x$ occurs on the branch; otherwise, it is called open. Closed and open branches are denoted respectively by $\otimes$ and $\odot$. A tableau closes if all branches of the tableau close. A tableau branch $\pi$ is \textit{complete} iff every rule that can be applied in $\pi$ has been applied.

A tableau proof of $\phi$ from a finite set of premises $T$ is a closed tableau beginning with a single branch, called the \textit{root} of the tableau, whose nodes consist of $\psi,+0$ for every premise $\psi \in T$ and $\phi,-0$, constructed by applying tableau rules from Figure \ref{Figure: Tableau rules of JRC}. We write $T \vdash_{\JRC} \phi$ if there is a tableau proof of $\phi$ from a set of premises $T$.

\begin{example}\label{example:JCR-tableau}
		We give a tableau proof for $\jbox{s} \phi \limpliesCounterRel  \jbox{s \ttsum t} \phi.$
		
		\vspace*{0.2cm}
		\Tree [.$1.~\jbox{s}\phi\limpliesCounterRel\jbox{s\ttsum t}\phi,-0$ [.$2.~0\relTableau{\jbox{s}\phi}1$ [.$3.~\jbox{s}\phi,+1$ [.$4.~\jbox{s+t}\phi,-1$ [.$5.~1\relTableau{s+t}2$ [.$6.~\phi,-2$ [.$7.~1\relTableau{s}2$ [.$8.~1\relTableau{t}2$ {$9.~\phi,+2$ \\ $\otimes$} ] ] ] ] ] ] ] ]
		\vspace*{0.2cm}
		
		\noindent
		Formula 1 is the root of the tableau. Formulas 2, 3 and 4 are from 1 by rule $(F\limpliesCounterRel_0)$, 5 and 6 are from 4  by rule $(F\jbox{})$,  7 and 8 are from 5 by $(\relTableau{+})$, and 9 is from 3 and 7 by rule $(T\jbox{})$.
\end{example}

Let us show the soundness and completeness of this tableau system. 

\begin{definition}
	Given a $\JRC$-model $\model= (\worlds, \worldsNormal, \relRelevant, \accrelFormula, \accrelTrem, \ast, \valuationR)$ and a tableau branch $\pi$, we say that $\model$ is \textit{faithful} to $\pi$ iff there exists a function $h: \mathbb{N} \to \worlds$ such that $h(0) \in \worldsNormal$, $h(i^\sharp) = (h(i))^\ast$ and the following hold:
	\begin{itemize}
		\item if $\phi, +x$ occurs in $\pi$, then $(\model, h(x)) \models \phi$,
		
		\item if $\phi, -x$ occurs in $\pi$, then $(\model, h(x)) \not \models \phi$,
		
		\item if $rxyz$ occurs in $\pi$, then $\relRelevant h(x)  h(y) h(z)$,
		
		\item if $x \relTableau{\phi} y$ occurs in $\pi$, then $h(x) \relModel{\phi} h(y)$,
		
		\item if $x \relTableau{t} y$ occurs in $\pi$, then $h(x) \relModel{t} h(y)$.
	\end{itemize}
\end{definition}

\begin{lemma}\label{lem: soundness lemma tableaux}
	Let $\pi$ be any branch of a tableau and let  $\model$ be a $\JRC$-model that is faithful to $\pi$. If a tableau rule is applied to $\pi$, then it produces at least one extension $\pi'$ such that $\model$ is faithful to $\pi'$.
\end{lemma}
\begin{proof}
	The proof is straightforward. Let us do only four cases in detail as an example. 
	
	Let $\model= (\worlds, \worldsNormal, \relRelevant, \accrelFormula, \accrelTrem, \ast, \valuationR)$ be a $\JRC$-model that is faithful to a branch $\pi$ of a tableau. Suppose that the rule 
	\[
	\alwaysNoLine\AXC{}
	\UIC{$$}
	\UIC{$\diagup\diagdown$}
	\UIC{$ \phi, -x~~~~~~~ \phi, +x$}
	\UIC{$~~~~~~~~~~~~~~~~~ x \relTableau{\phi} x$}
	\UIC{}
	\DP
	\]
	is applied to $\pi$. Then, we have either $(\model, h(x)) \not \models \phi$ or $(\model, h(x)) \models \phi$. In the latter case, by condition 2 of Definition \ref{def: relational models conditions JRC}, we get $h(x) \accrel_{\phi} h(x)$. Thus, either $\model$ is faithful to $\pi$ extended by $\phi, -x$ or it is faithful to $\pi$ extended by $\phi, +x$ and $x \relTableau{\phi} x$.
	
	Suppose that the rule 
	\[
	\alwaysNoLine
	\AXC{}
	\UIC{$\phi \limpliesCounterRel \psi, -0$}
	\UIC{$|$}
	\UIC{$0 \relTableau{\phi} j$}
	\UIC{$ \phi, +j$}
	\UIC{$\psi, -j$}
	\UIC{}
	\DP
	\]
	is applied to $\pi$, where $j$ is new. Call the new branch $\pi'$. Thus, $(\model, h(0)) \not \models \phi \limpliesCounterRel \psi$, and hence $\relModel{\phi} (h(0)) \not \subseteq \truthsetModel{\psi}$. Therefore, there exists $v \in \relModel{\phi} (h(0))$ such that $v \not \in \truthsetModel{\psi}$. Since $h(0) \in \worldsNormal$, by condition 1 of Definition \ref{def: relational models conditions JRC}, we get $v \in \truthsetModel{\phi}$. Let $h'$ be the same function as $h$ except that $h'(j) = v$. Then, since $j$ was new, $h'$ is faithful to $\pi$. Therefore we have $h'(0) \relModel{\phi} h'(j)$, $(\model, h'(j)) \models \phi$, and $(\model, h'(j)) \not \models \psi$. Thus, $\model$ is faithful to $\pi'$.
	
	Suppose that the rule 
	\[
	\alwaysNoLine\AXC{}
	\UIC{$\jbox{t} \phi, +x$}
	\UIC{$x \relTableau{t} y$}
	\UIC{$|$}
	\UIC{$ \phi, +y$}\DP
	\]
	is applied to $\pi$  and a new branch $\pi'$ is produced. Thus, $(\model, h(x)) \models \jbox{t} \phi$ and $h(x) \relModel{t} h(y)$. Then, it follows that $(\model, h(y)) \models \phi$. Thus, $\model$ is faithful to $\pi'$.
	
	Suppose that the rule 
	\[
	\alwaysNoLine
	\AXC{}
	\UIC{$\jbox{t} \phi, -x$}
	\UIC{$|$}
	\UIC{$x \relTableau{t} j$}
	\UIC{$\phi, -j$}
	\DP
	\]
	is applied to $\pi$, where $j$ is new, and a new branch $\pi'$ is produced. Thus, $(\model, h(x)) \not \models \jbox{t} \phi$, and hence $\relModel{t} (h(x)) \not \subseteq \truthsetModel{\phi}$. Therefore, there exists $v \in \relModel{t} (h(x))$ such that $v \not \in \truthsetModel{\phi}$.  Let $h'$ be the same function as $h$ except that $h'(j) = v$. Then, since $j$ was new, $h'$ is faithful to $\pi$. Therefore we have $h'(x) \relModel{t} h'(j)$ and   $(\model, h'(j)) \not \models \phi$. Thus, $\model$ is faithful to $\pi'$.
	\qed
\end{proof}

Now soundness is simple to establish.

\begin{theorem}[Soundness]\label{thm:Soundness tableaux}
	Given a finite set of premises $T$, if $T \vdash_{\JRC} \phi$, then $T \models_{\JRC} \phi$. 
\end{theorem}
\begin{proof}
	If $T \not \models_{\JRC} \phi$, then there is a $\JRC$-model $\model= (\worlds, \worldsNormal, \relRelevant, \accrelFormula, \accrelTrem, \ast, \valuationR)$ and a normal world $w \in \worldsNormal$ such that $(\model, w) \entails \psi$ for all $\psi \in T$ and $(\model, w) \not \entails \phi$.  Thus, by Lemma~\ref{lem: soundness lemma tableaux}, there is no tableau proof of $\phi$ from $T$, since otherwise $\model$ would be faithful to closed branches, which is impossible. Therefore, $T \not \vdash_{\JRC} \phi$. \qed
\end{proof}

The proof of completeness is similar to that given in~\cite{Priest2008}, and so the details are omitted here.

\begin{definition}
	Let $\pi$ be an open complete branch of a tableau. The model $\model_\pi$ induced by $\pi$ is defined as follows:
	\begin{itemize}
		\item $\worlds := \{ w_x \mid x = i $ or $i^\sharp$, for some $i \in \mathbb{N}$, and $x$ occurs on $\pi \}$ and $\worldsNormal := \lset{ w_0 }$.
		
		\item $(w_x)^\ast := w_{\bar{x}}$.
		
		\item $\relRelevant w_x w_y w_z$ iff $r xyz$ occurs on $\pi$.
		
		\item if $\phi$ occurs as the antecedent of a conditional $\limpliesCounterRel$ at a node
		of $\pi$, then $w_x \relModel{\phi} w_y$ iff $x \relTableau{\phi} y$ is on $\pi$;
		otherwise, $w_x \relModel{\phi} w_y$ iff $(\model_\pi, w_y) \models \phi$.
		
		\item $w_x \relModel{t} w_y$ iff $x \relTableau{t} y$ occurs on $\pi$.
		
		\item $\valuationR(p) = \{ w_x \mid p, +x$ occurs on $\pi \}$, for $p \in \Prop$.
	\end{itemize}
\end{definition}

\begin{lemma}\label{lem: completeness lemma tableaux}
	Let $\pi$ be any open complete branch of a tableau and let  $\model_\pi$ be the model induced by $\pi$. Then,
	\begin{enumerate}
		\item If $\phi, +x$ occurs in $\pi$, then $(\model_\pi, w_x) \models \phi$.
		\item If $\phi, -x$ occurs in $\pi$, then $(\model_\pi, w_x) \not \models \phi$.
	\end{enumerate}
\end{lemma}
\begin{proof}
	The proof is by induction on the complexity of $\phi$. If $\phi$ is an atomic proposition, then the result is true by the definition of $\model_\pi$. The proof for the case when $\phi$ is either $\lnegRelevant \psi$ or $\psi \wedge \chi$ or $\psi \limpliesRelevant \chi$ is standard (cf. \cite{Priest2008}). 
	
	Suppose that $\phi = \psi \limpliesCounterRel \chi$ and $\psi \limpliesCounterRel \chi, +x$ occurs in $\pi$. We have two cases:
	\begin{description}
		\item[$(i)$] For no $y$, $x \relTableau{\psi} y$ occurs in $\pi$. Then, $\relModel{\psi}  (w_x)= \emptyset$, and hence $\relModel{\psi} (w_x) \subseteq \truthset{\chi}{{\model_\pi}}$. Thus, $(\model_\pi, w_x) \models \psi \limpliesCounterRel \chi$.
		
		\item[$(ii)$] For some $y$, $x \relTableau{\psi} y$ occurs on $\pi$. Let $w_y \in \relModel{\psi} (w_x)$ be an arbitrary state. Since $\psi$ occurs as the antecedent of a conditional $\limpliesCounterRel$ at a node
		of $\pi$, $x \relTableau{\psi} y$ occurs on $\pi$. Since $\pi$ is a complete branch, $\chi, +y$ occurs on $\pi$. By the induction hypothesis, $(\model_\pi, w_y) \models \chi$. Hence, $\relModel{\psi} (w_x) \subseteq \truthset{\chi}{{\model_\pi}}$, and thus $(\model_\pi, w_x) \models \psi \limpliesCounterRel \chi$.
	\end{description}
	
	Now suppose that $\psi \limpliesCounterRel \chi, -x$ occurs on $\pi$. Since $\pi$ is a complete branch, $x \relTableau{\psi} j$ and $\chi, -j$ occurs in $\pi$, for some  $j \in \mathbb{N}$. Since $\psi$ occurs as the antecedent of a conditional $\limpliesCounterRel$ at a node of $\pi$, $w_x \relModel{\psi} w_j$. On the other hand, by the induction hypothesis, $(\model_\pi, w_j) \not \models \chi$. Hence, $\relModel{\psi} (w_x) \not \subseteq \truthset{\chi}{{\model_\pi}}$, and thus $(\model_\pi, w_x) \not \models \psi \limpliesCounterRel \chi$.
	
	Next suppose that $\phi = \jbox{t} \psi$ and $\jbox{t} \psi, +x$ occurs on $\pi$. Let $w_y \in \relModel{t} (w_x)$ be an arbitrary state. Hence, $x \relTableau{t} y$ occurs on $\pi$. Since $\pi$ is a complete branch, $\psi, +y$ occurs on $\pi$. By the induction hypothesis, $(\model_\pi, w_y) \models \psi$. Hence, $\relModel{t} (w_x) \subseteq \truthset{\psi}{{\model_\pi}}$, and thus $(\model_\pi, w_x) \models \jbox{t} \psi$.
	
	Now suppose that $\jbox{t} \psi, -x$ occurs on $\pi$. Since $\pi$ is a complete branch, $x \relTableau{t} j$ and $\psi, -j$ occurs in $\pi$, for some  $j$. Hence, $w_x \relModel{t} w_j$ and, by the induction hypothesis, $(\model_\pi, w_j) \not \models \psi$. Hence, $\relModel{t} (w_x) \not \subseteq \truthset{\psi}{{\model_\pi}}$, and thus $(\model_\pi, w_x) \not \models \jbox{t} \psi$.
	\qed
\end{proof}

\begin{lemma}\label{lem: induced model is JCR-model}
	Let $\pi$ be any open complete branch of a tableau and let  $\model_\pi$ be the model induced by $\pi$. Then, $\model_\pi$ is a $\JRC$-model.
\end{lemma}
\begin{proof}
	First note that for any $w_x \in \worlds$ we have that $(w_x)^{**} = (w_{\bar{x}})^* = w_{\bar{\bar{x}}} = w_x$. The proof of the normality condition is similar to that given in \cite{Priest2008}[page 212]. Next, we show that $\model_\pi$ satisfies conditions 1--3 of Definition \ref{def: relational models conditions JRC}. 
	\begin{enumerate}
		\item  Let $w_y \in \relModel{\phi} (w_0)$ be an arbitrary state. We have to show that $(\model_\pi, w_y) \models \phi$. We distinguish the following possibilities:
		\begin{description}
			\item[$(i)$] $\phi$ occurs as the antecedent of a conditional $\limpliesCounterRel$ at a node of $\pi$. Hence, by the definition of $\model_\pi$,  $0 \relTableau{\phi} y$ occurs on $\pi$. Then, we have two cases. 
			\begin{description}
				\item[$(a)$] $0 \relTableau{\phi} y$ is produced by the rule
				\[
				\alwaysNoLine
				\AXC{}
				\UIC{$\phi \limpliesCounterRel \psi, -0$}
				\UIC{$|$}
				\UIC{$0 \relTableau{\phi} y$}
				\UIC{$ \phi, +y$}
				\UIC{$\psi, -y$}
				\UIC{}
				\DP
				\]
				where $y \in \N$. In this case, since $ \phi, +y$ occurs on $\pi$, by Lemma \ref{lem: completeness lemma tableaux}, we get $(\model_\pi, w_y) \models \phi$, which is what we wished to show. 
				
				\item[$(b)$]  $0 \relTableau{\phi} y$ is produced by the rule
				\[
				\alwaysNoLine\AXC{}
				\UIC{$$}
				\UIC{$\diagup\diagdown$}
				\UIC{$ \phi, -0~~~~~~~ \phi, +0$}
				\UIC{$~~~~~~~~~~~~~~~~~ 0 \relTableau{\phi} 0$}
				\UIC{}
				\DP
				\]
				where $y = 0$ and $\pi$ is the right branch. In this case, since $ \phi, +0$ occurs in $\pi$, by Lemma \ref{lem: completeness lemma tableaux}, we get $(\model_\pi, w_0) \models \phi$.
				
			\end{description}

			\item[$(ii)$] $\phi$ does not occur as the antecedent of a conditional $\limpliesCounterRel$ on $\pi$. Then, by the definition of $\model_\pi$, $w_0 \relModel{\phi} w_y$ iff $(\model_\pi, w_y) \models \phi$.
			
		\end{description}
		Thus, in all cases we get $(\model_\pi, w_y) \models \phi$. Since $w_y$ was an arbitrary state in $\relModel{\phi} (w_0)$, we have that $\relModel{\phi} (w_0) \subseteq \truthset{\phi}{{\model_\pi}}$.
		
		\item Let $w_x \in \worlds$ and $(\model_\pi, w_x) \models \phi$. We have to show that $w_x \in \relModel{\phi} (w_x)$. We distinguish again the following possibilities:
		\begin{description}
			\item[$(i)$] $\phi$ occurs as the antecedent of a conditional $\limpliesCounterRel$ at a node of $\pi$. Since $w_x \in \worlds$ and $\pi$ is a complete branch, $x$ occurs on $\pi$ and the following rule $(Cut_r)$ is applied on $\pi$
			\[
			\alwaysNoLine\AXC{}
			\UIC{$$}
			\UIC{$\diagup\diagdown$}
			\UIC{$ \phi, -x~~~~~~~ \phi, +x$}
			\UIC{$~~~~~~~~~~~~~~~~~ x \relTableau{\phi} x$}
			\UIC{}
			\DP
			\]
			Note that, since $(\model_\pi, w_x) \models \phi$, the right branch is $\pi$. Thus, $x \relTableau{\phi} x$ occurs on $\pi$, and therefore $w_x \relModel{\phi} w_x$.

			\item[$(ii)$] $\phi$ does not occur as the antecedent of a conditional $\limpliesCounterRel$ on $\pi$. Then, by the definition of $\model_\pi$,  from $(\model_\pi, w_x) \models \phi$ it follows that $w_x \relModel{\phi} w_x$.
			
		\end{description}
		
		\item Given $s,t \in \TermsJCR$, let $w_x \relModel{s+t} w_y$. Then, by the definition of $\model_\pi$, $x \relTableau{s+t} y$ occurs on $\pi$. Since $\pi$ is a complete branch, by applying the rule $(\relTableau{+})$, both $x \relTableau{s} y$ and $x \relTableau{t} y$ occurs on $\pi$. Thus, $x \relModel{s} y$ and $x \relModel{t} y$.

	\end{enumerate}
	\qed
\end{proof}

\begin{theorem}[Completeness]\label{thm:completeness tableaux}
	Given a finite set of premises $T$, if  $T \models_{\JRC} \phi$, then $T \vdash_{\JRC} \phi$. 
\end{theorem}
\begin{proof}
	If $T \not \vdash_{\JRC} \phi$, then there is an open complete branch, say $\pi$ in the tableau of $\phi$ from $T$. Let $\model_\pi$ be the model induced by $\pi$. By Lemma \ref{lem: induced model is JCR-model}, $\model_\pi$ is a $\JRC$-model. By Lemma \ref{lem: completeness lemma tableaux},  $(\model_\pi, w) \entails \psi$ for all $\psi \in T$ and $(\model_\pi, w) \entails \phi$. Thus, $T \not \models_{\JRC} \phi$. 
\end{proof}

\paragraph{\textbf{Counterpossibles in $\JRC$}.}

Using tableaux it is easy to show that neither $(p \wedge \lnegRelevant p) \limpliesCounterRel \phi$ nor $\phi \limpliesCounterRel (p \vee \lnegRelevant p)$, for $p \in \Prop$, is valid in $\JRC$, and thus counterpossibles are not valid.

\begin{lemma}\label{lem: counterpossible in JCR}
	For some $\phi \in \FormulaeJCR$ and $p \in \Prop$, we have:
	\begin{enumerate}
		\item $\not \entails_{\JRC} (p \wedge \lnegRelevant p) \limpliesCounterRel \phi$.
		
		\item $\not \entails_{\JRC} \phi \limpliesCounterRel (p \vee \lnegRelevant p)$.
	\end{enumerate}
\end{lemma}
\begin{proof}
	By constructing a tableau for the formula $(p \wedge \lnegRelevant p) \limpliesCounterRel q$, we show that it is not valid in $\JRC$:
	\[
	\alwaysNoLine
	\AXC{$(p \wedge \lnegRelevant p) \limpliesCounterRel q, -0$}
	\UIC{$|$}
	\UIC{$0 \relTableau{p \wedge \lnegRelevant p} 1$}
	\UIC{$p \wedge \lnegRelevant p, +1$}
	\UIC{$q, -1$}
	\UIC{$p, +1$}
	\UIC{$\lnegRelevant p, +1$}
	\UIC{$p, - 1^\sharp$}
	\UIC{$r000$}
	\UIC{$r011$}
	\UIC{$r01^\sharp 1^\sharp$}
	\UIC{$\odot$}
	\DP
	\]
	This branch is open, and it cannot be closed even by more applications of $(F\limpliesCounterRel_0)$ or the rule $(Cut_r)$. In fact, applying  the rule $(Cut_r)$ will produce open branches containing either $p \wedge \lnegRelevant p, -0$, or $p \wedge \lnegRelevant p, +0$, or $1 \relTableau{p \wedge \lnegRelevant p} 1$ or $p \wedge \lnegRelevant p, -1^\sharp$.  And applying more applications of $(F\limpliesCounterRel_0)$ produces nodes of the form $0 \relTableau{p \wedge \lnegRelevant p} j$ and $p \wedge \lnegRelevant p, +j$, for $j \in \mathbb{N}$, that can be treated similar to above branch with no closure. The following Routley relational model $\model = (\worlds, \worldsNormal, \relRelevant, \accrelFormula, \accrelTrem, *,  \valuationR)$ is constructed from one of these open complete branches in which $0^\sharp$ does not occur:
	\begin{itemize}
		\item $\worlds =  \{ w_0, w_1, w_1^*\}$, 
		
		\item $\worldsNormal = \{ w_0 \}$,
		
		\item $w_0^* = w_0$,
		
		\item $\relRelevant = \lset{ (w_0,w_0,w_0), (w_0,w_1,w_1), (w_0,w_1^*,w_1^*) }$
		
		\item $\accrel_{\phi}(w_0) = \accrel_{\phi}(w_1) =  \accrel_{\phi}(w_1^*) =\truthsetModel{\phi}$, for all $\phi \in \FormulaeJCR$, 
		
		\item $\accrel_{t}(w_0) = \accrel_{t}(w_1) = \accrel_{t}(w_1^*) = \emptyset$, for all $t \in \TermsJCR  $,
		
		\item $\valuationR(p) = \lset{w_1}$ and $\valuationR(q) = \emptyset$. $\valuationR$ is defined arbitrary for other atomic propositions.
	\end{itemize}
	Clearly, $\model$ is a $\JRC$-model. Since $\accrel_{p \wedge \lnegRelevant p}(w_0) = \truthset{p \wedge \lnegRelevant p}{\model} = \lset{w_1}$ and $\truthset{q}{\model} = \emptyset$, we get $(\model, w_0) \not \models (p \wedge \lnegRelevant p) \limpliesCounterRel q$, as desired. 
	
	Now consider the model $\model' = (\worlds, \worldsNormal, \relRelevant,\accrelFormula, \accrelTrem, *,  \valuationR')$ which is defined similar to $\model$ except that the valuation $\valuationR'$ is defined as follows:
	\begin{itemize}
		\item $\valuationR'(p) = \lset{w_1^*}$ and $\valuationR'(q) = \lset{w_1}$. $\valuationR'$ is defined arbitrarily for other atomic propositions.  
	\end{itemize}
	Since $\relModel{q}(w_0)= \truthsetModel{q} = \lset{w_1} \not \subseteq \truthsetModel{p \vee \lnegRelevant p} = \lset{w_0, w_1^*}$, we get $(\model, w_0) \not \models q \limpliesCounterRel (p \vee \lnegRelevant p)$.
	\qed
\end{proof}

\section{Justified true belief, again}
\label{sec: JTB in JCR}

In this section, we study the definitions of knowledge in the logic $\JRC$. Both the justified true belief and Nozick's definition of knowledge can be expressed in the language of $\JRC$. The justified true belief of $\phi$ by justification $t$ can be expressed by $\jbox{t} \phi$, and
Nozick's definition of knowledge can be formalized by the following formula 
\[
\phi \wedge \jbox{t} \phi \wedge (\lnegRelevant \phi \limpliesCounterRel \lnegRelevant \jbox{t} \phi) \wedge (\phi \limpliesCounterRel \jbox{t} \phi).
\]
The Aumann knowledge, formalized by the universal modality $\lknows$, can be added to the language of $\JRC$ similar to that defined in Section \ref{sec: Nozick's analysis of knowledge}. By formalizing a Gettier-type example, Chisholm's example ``the sheep in the field'' from \cite{Chisholm-1966}, we show that neither the justified true belief nor Nozick's definition of knowledge is Aumann knowledge. 

\begin{example}[Chisholm's example]
	``Imagine that you are standing outside a field. You see, within it, what looks exactly like a sheep. What belief instantly occurs to you? Among the many that could have done so, it happens to be the belief that there is a sheep in the field. And in fact you are right, because there is a sheep behind the hill in the middle of the field. You cannot see that sheep, though, and you have no direct evidence of its existence. Moreover, what you are seeing is a dog, disguised as a sheep. Hence, you have a well justified true belief that there is a sheep in the field. But is that belief knowledge?'' (Quoted  from \cite{Hetherington-IEP})
	
	Let us formalize this example in $\JRC$. Consider the following dictionary, where $p,q \in \Prop$ and $\justVarOne \in \VTerms$:
	\begin{align*}
		p &:=  \text{ There is a sheep in the field.} \\
		q &:=  \text{ The animal I'm seeing is a sheep.} \\
		\justVarOne &:=  \text{ My observation of the animal.} 
	\end{align*}

	\begin{figure}
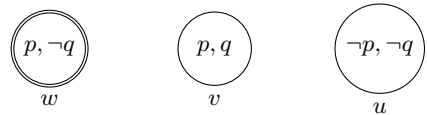

		\begin{center}
			\begin{tabular}{c}
				\begin{mytikz}
					
					\node[w,double,label={below:$w$}] (w) {$p, \neg q$};
					
					\node[w,right of=w,label={below:$v$}] (v) {$p, q$};
					
					\node[w,right of=v,label={below:$u$}] (u) {$\neg p, \neg q$};

				\end{mytikz}
			\end{tabular}
		\end{center}
		\caption{Model $\mathcal{N}$ for the sheep in the field example. $w$ is the actual state.}
	\end{figure}
	Define the Routley relational model $\mathcal{N} = (\worlds, \worldsNormal, \relRelevant, \accrelFormula, \accrelTrem, *,  \valuationR)$ as follows:
	\begin{itemize}
		\item $\worlds = \worldsNormal = \{ w, v, u\}$, 
		
		
		\item $w^* = w$, $v^* = v$,  and $u^* = u$,
		
		\item $\relRelevant = \worlds \times \worlds \times \worlds$.
		
		\item $\accrel_{\phi}(w) = \accrel_{\phi}(v) =  \accrel_{\phi}(u) =\truthsetModel{\phi}$, for all $\phi \in \FormulaeJCR$,

		\item $\accrel_{t}(w) = \lset{w}$, $\accrel_{t}(v) = \lset{v}$ and $\accrel_{t}(u) = \lset{u}$, for all $t \in \TermsJCR  $,

		\item $\valuationR(p) = \lset{w, v}$, $\valuationR(q) = \lset{v}$. $\valuationR$ is defined arbitrarily for other atomic propositions.
	\end{itemize}
	It is easy to show that $\mathcal{N}$ is an $\JRC$-model. We observe that $\mathcal{N} \models q \limpliesCounterRel p$, which means that the conditional ``if the animal I'm seeing were a sheep, then there would be a sheep in the field'' is true, as expected.
	
	Since $\accrel_{\justVarOne}(w) = \lset{w} \subseteq \truthset{p}{{\mathcal{N}}} = \lset{w,v}$, we get $(\mathcal{N}, w) \models \jbox{\justVarOne} p$. In fact, $\truthset{\jbox{\justVarOne} p}{{\mathcal{N}}} = \lset{w, v}$.
	
	On the other hand, since $\accrel_{\lnegRelevant p}(w) = \truthset{\lnegRelevant  p}{{\mathcal{N}}} = \lset{u}$ and $\truthset{\lnegRelevant \jbox{\justVarOne} p}{{\mathcal{N}}} = \lset{u}$, we have that  $(\mathcal{N},w)  \entails \lnegRelevant p \limpliesCounterRel \lnegRelevant \jbox{\justVarOne} p$. In addition, since  $\accrel_{p}(w) = \truthset{p}{{\mathcal{N}}} = \lset{w,v}   \subseteq \truthset{\jbox{\justVarOne} p}{{\mathcal{N}}}$, we get $(\mathcal{N},w)  \entails  p \limpliesCounterRel  \jbox{\justVarOne} p$.

	We have also that $(\mathcal{N},u) \not \entails p$ and thus $(\mathcal{N},w) \not \entails \lknows p$. Therefore, we showed that 
	\[
	(\mathcal{N}, w) \models p \wedge \jbox{\justVarOne} p \wedge (\lnegRelevant \phi \limpliesCounterRel \lnegRelevant \jbox{t} \phi) \wedge (\phi \limpliesCounterRel \jbox{t} \phi)\wedge \lnegRelevant \lknows p.
	\]
	These facts show that neither the justified true belief nor Nozick's definition of knowledge is  Aumann knowledge and also show that neither the sensitivity condition nor the adherence condition (formalized via the conditional $\limpliesCounterRel$) is sufficient for Aumann knowledge.
\end{example}


\bibliography{library}
\end{document}